\documentclass[12pt]{article}
\usepackage{amssymb}
\usepackage{latexsym}
\usepackage{amsfonts}
\usepackage{amsmath}
\newtheorem{theorem}{Theorem}[section]
\newtheorem{lemma}[theorem]{Lemma}
\newtheorem{proposition}[theorem]{Proposition}

\newtheorem{corollary}[theorem]{Corollary}

\numberwithin{equation}{section}

\usepackage[dvips]{color}
\usepackage{xcolor}
\usepackage{comment}

\def\P{\mathbb{P} }
\def\E{\mathbb{E} }
\def\R{\mathbb{R} }
\def\N{\mathbb{N} }

\def\tl{\tilde}
\def\d{{\rm d}}

\def\bP{{\bf P} }
\def\bQ{{\bf Q} }
\def\mP{{\mbox{P}} }
\def\mE{{\mbox{E}} }

\def\bE{{\bf E} }

\def\cL{{\cal L} }

\def\cM{{\cal M}}

\begin{document}

\allowdisplaybreaks

\title{\bf Lower
deviation for the supremum of the support of super-Brownian motion
\footnote{The research of this project is supported by the National Key R\&D Program of China (No. 2020YFA0712900).}
     }
\author{{\bf Yan-Xia Ren}\thanks{The research of this author is supported by   NSFC (Grant Nos. 11731009, 12071011 and  12231002) and LMEQF.}
\quad
{\bf Renming Song}\thanks{Research supported in part by a grant from the Simons
Foundation (\#960480, Renming Song).}
\quad
{\bf Rui Zhang}\thanks{The research of this  author is  supported by  NSFC (Grant Nos. 11601354 and 12271374), Beijing Municipal Natural Science Foundation(Grant No. 1202004), and Academy for Multidisciplinary Studies, Capital Normal University.}
}

\date{}
\maketitle

\begin{abstract}
We study the asymptotic behavior
of the supremum $M_t$ of the support of a supercritical super-Brownian motion.
In our recent paper (Stoch. Proc. Appl. \textbf{137} (2021), 1--34), we showed that,  under some conditions,
$M_t-m(t)$ converges
in distribution to a randomly shifted Gumbel random variable,
where $m(t)=c_0t-c_1\log t$. In the same paper, we also studied the upper large deviation of $M_t$, i.e., the asymptotic behavior of $\P(M_t>\delta c_0t) $ for $\delta\ge 1$.
In this paper, we study the lower large deviation of $M_t$, i.e., the asymptotic behavior of $\P(M_t\le \delta c_0t|\mathcal{S}) $ for $\delta<1$, where $\mathcal{S}$ is the survival event.

\end{abstract}

\medskip


	\begin{doublespace}
	\section{Introduction}
	\subsection{Super-Brownian motion}
	Let $\psi$ be a function of the form
	$$\psi(\lambda)=-\alpha\lambda+\beta\lambda^2+\int_0^\infty \Big(e^{-\lambda y}-1+\lambda y\Big)n(\d y),\quad \lambda\ge 0,$$
	where
	$\alpha\in\R$,
	$\beta \ge 0$ and $n$ is a $\sigma$-finite measure satisfying
	$$
	\int_0^\infty (y^2\wedge y) n(\d y)<\infty.
	$$
	$\psi$ is called a branching mechanism.
	We will always assume that $\lim_{\lambda\to\infty}\psi(\lambda)=\infty$.
Let $\{B_t, t\geq 0; \mbox{P}_x\}$ be a standard Brownian motion starting from $x\in\R$,
and let $\mE_x$ be the corresponding expectation.
 We write  $\mP=\mP_{0}$
  and  $\mE=\mE_0$.
	In this paper we will consider a super-Brownian motion $X$ on $\R$ with branching mechanism $\psi$.

	Let $\mathcal{B}^+(\R)$ (resp. $\mathcal{B}^+_b(\R)$) be the space of non-negative
	(resp. bounded non-negative)
	Borel functions on $\R$,
	and let ${\cal M}_F(\R)$
	be the space of finite measures on $\R$, equipped with the topology of weak convergence. A super-Brownian motion
$X=\{X_t,t\geq 0\}$
with branching mechanism $\psi$ is a  Markov process taking values in ${\cal M}_F(\R)$.
	For any $\mu \in \mathcal{M}_F(\R)$, we denote the
	law of $X$ with initial configuration $\mu$
	by $\P_\mu$, and the corresponding expectation by $\E_\mu$.
We write $\P=\P_{\delta_0}$
and $\E=\E_{\delta_0}$.
	As usual, we use the notation
	$\langle f,\mu\rangle:=\int_{\R} f(x)\mu(dx)$
	and $\|\mu\|:=\langle 1,\mu\rangle$. Then for all
	$f\in \mathcal{B}^+_b(\R)$ and
$\mu \in \mathcal{M}_F(\R)$,
	\begin{equation}\label{V}
		-\log \E_\mu\left(e^{-\langle f,X_t\rangle}\right)=
		\langle V_f(t, \cdot),\mu\rangle, \qquad t\geq 0,
	\end{equation}
	where $V_f(t, x)$ is the unique positive solution to the equation
	\begin{equation}\label{eqt-u}
		V_f(t, x)+\mE_x\int_0^t\psi(V_f(t-s, B_s))\d s=\mE_x f(B_t),\qquad t\geq 0.
	\end{equation}
The existence of such superprocesses is well-known, see, for instance,
\cite{Dawson}, \cite{E.B.} or \cite{Li11}.

	It is well known that $\|X_t\|$ is a continuous state branching process with branching mechanism $\psi$ and that
	$$\P(\lim_{t\to\infty}\|X_t\|=0)=e^{-\lambda^*},$$
	where $\lambda^*\in[0,\infty)$ is the largest root of the equation $\psi(\lambda)=0$.
        It is known that $\lambda^*>0$ if and only if $\alpha=-\psi'(0+)>0$.
	$X$ is called a supercritical (critical, subcritical) super-Brownian motion if $\alpha > 0$ ($= 0,< 0$).
	In this paper, we only deal with the supercritical case, that is, we assume $\alpha > 0$.
	Let $M_t$ be the supremum of the support of $X_t$.
More precisely,
we define the rightmost point  $M(\mu)$ of
$\mu\in\cM_{F}(\R)$ by
$M(\mu):=\sup\{x: \mu(x,\infty)>0\}$.
Here we use the convention that $\sup\emptyset=-\infty.$
Then $M_t$ is simply $M(X_t)$.
Recently,  in \cite{SRZ}, we studied the
asymptotic behavior of $M_t$ under the following two assumptions:
	\begin{itemize}
		\item[] {\bf{(H1)}} There exists $\gamma>0$ such that
		
		\begin{equation*}\label{cond-log}
			\int_1^\infty y(\log y)^{2+\gamma}n(\d y)<\infty.
		\end{equation*}
		\item[]{\bf (H2)}
		There exist $\vartheta\in(0,1]$ and $a>0,b>0$ such that
		\begin{equation*}\label{H2}
			\psi(\lambda)\ge -a\lambda+b\lambda^{1+\vartheta}, \quad \lambda >0.
		\end{equation*}
	\end{itemize}
It is clear that if  $\beta>0$ or $n(\d y)\ge y^{-1-\vartheta}\,\d y$,  then (H2) holds.
	Condition (H2) implies that the following Grey condition holds:
\begin{equation}\label{Grey}
	\int^\infty \frac{1}{\psi(\lambda)}\ d\lambda<\infty.
\end{equation}
It is well known that under the above Grey condition,
$\lim_{t\to\infty}\P_{\mu}(\|X_t\|=0)=e^{-\lambda^*\|\mu\|}.$
Denote  $\mathcal{S}:=\{\forall t\ge 0, \|X_t\|>0\}$. It is clear that $\P(\mathcal{S})\in(0,1)$.
Define, for $t\ge0$,
$$
D_t:=\langle (\sqrt{2\alpha}t-\cdot)e^{-\sqrt{2\alpha}(\sqrt{2\alpha}t-\cdot)},X_t \rangle.
$$
It has been proven in \cite{KLSR}
that $\{D_t, t\geq 0\}$ is a martingale,
which is called the derivative martingale of the super-Brownian motion $X_t$,
and that $D_t$ has an almost sure non-negative limit $D_\infty$ as $t\to\infty$.
	Assumption (H2) also implies that
	\begin{eqnarray}\label{con:psi}
		\int^\infty\frac{1}{\sqrt{\int_{\lambda^*}^\xi \psi(u)\,du}}\,\d \xi<\infty.
	\end{eqnarray}
	Under (H1) and \eqref{con:psi},   $D_\infty$ is non-degenerate and
	\begin{equation}\label{limit-as}
		\frac{M_t}{t}\to \sqrt{2\alpha},\quad \P\mbox{-a.s. on } \mathcal{S},
	\end{equation}	
	see \cite[Theorem 2.4 and Corollary 3.2 ]{KLSR}.
	
		For any $f\in \mathcal{B}^+(\R)$, put
	\begin{equation}\label{def:u}
		u_f(t,x):=-\log \E \left( e^{-\int_\R f(y-x) X_t(dy)}; M_t\le x\right),
	\end{equation}
Note that $u_f$ only depends on the value of $f$ on $(-\infty,0]$.
Let $\mathcal{H}$ be the space of all the nonnegative bounded functions $f$ on
$(-\infty,0]$ satisfying
\begin{equation}\label{initial-cond1'}
	\int_0^\infty y e^{\sqrt{2\alpha}y}f(-y)\,dy<\infty.
\end{equation}
It has been proved in \cite[Theorem 1.3]{SRZ} that
under  (H1)-(H2),
for any $f\in \mathcal{H}$,
we have that
	\begin{equation}\label{M-w}
		\lim_{t\to\infty}u_f(t,m(t)+x)=w_f(x),
	\end{equation}
	where
	\begin{equation}\label{def-m_t}
		m_t=\sqrt{2\alpha}t-\frac{3}{2\sqrt{2\alpha}}\log t,
	\end{equation}
and
$w_f$ is a traveling wave solution of the F-KPP equation,
that is, a solution of
 $$
  \frac{1}{2}w_{xx}+\sqrt{2\alpha}w_x-\psi(w)=0.
 $$
 Moreover,
 $w_f$ is given by
	$w_f(x)=-\log \E\left[\exp\{-\tilde{C}(f)D_\infty e^{-\sqrt{2\alpha}x}\}\right]$,
 with
$$\tl{C}(f):=\lim_{r\to\infty}\sqrt{\frac{2}{\pi}}\int_0^\infty u_{f}(r,\sqrt{2\alpha}r+y)ye^{\sqrt{2\alpha}y}\,dy\in(0,\infty).$$

In the remainder of this paper,
we write $u(t,x)$ and $w(x)$ for $u_f(t,x)$ and $w_f(x)$
respectively when $f\equiv 0$.
	
\subsection{Main results}
	
   In \cite[Theorem 1.2]{SRZ},
   we proved the following upper large deviation results for $M_t$ under conditions (H1)-(H2):
	\begin{itemize}
		\item[(1)]For $\delta>1$,
		$$\lim_{t\to\infty}\sqrt{t}e^{\alpha(\delta^2-1)}\P(M_t>\sqrt{2\alpha}\delta t)\in(0,\infty);$$
		\item[(2)]
		$$\lim_{t\to\infty}\frac{t^{3/2}}{\frac{3}{2\sqrt{2\alpha}}\log t}\P(M_t>\sqrt{2\alpha}t)\in(0,\infty).$$
	\end{itemize}
	However, using the methods in \cite{SRZ}, we
could not get the asymptotic behavior of the lower large deviation probability $\P(M_t\le \sqrt{2\alpha}\delta t|\mathcal{S})$ for $\delta<1$.
The purpose of this paper is to study the asymptotic behavior of the lower large deviation probability.
To accomplish this, we use the skeleton decomposition of super-Brownian motion and adapt some ideas from \cite{CHM}
used in the study of lower deviations of  the maximum of branching Brownian motion.
	
	  For branching Brownian motion, the aysmptotic behavior of
	  the maximal position, also denoted by $M_t$,  of the particles alive at time $t$
	  has been intensively studied.
	 	   To simplify notation, we consider a standard binary branching Brownian motion in $\R$, i.e.,
	  the lifetime of a particle is an exponential random variable with parameter 1  and when it dies,  it gives birth to $2$ children at the position of its death.
	 	  Bramson  proved in \cite{Bramson78} that
 $P(M_t-m(t)\le x)\to 1-w(x)$
	  as $t\to\infty$,
	  	  where $m(t)=\sqrt{2}t-\frac{3}{2\sqrt{2}}\log t$
	  and  $w(x)$ is a traveling wave solution. For the large deviation of $M_t$, \cite{ Chauvin88,Chauvin} studied the convergence rate of
 $P(M_t>\sqrt{2}\delta t)$ for $\delta \ge 1$.
 Recently,  Derrida and Shi \cite{DS1, DS}
 studied the lower large deviation of $M_t$, i.e, the asymptotic behavior of $\frac1{t}\log P(M_t\le \sqrt{2}\delta t)$ for $\delta <1$, and found that the rate function has a phase transition at $1-\sqrt{2}$.
 In \cite{CHM}, Chen, He and Mallein studied the limiting property of
 $P(M_t\le \sqrt{2}\delta t)$ for $\delta<1.$
 For more results on extremal processes of branching Brownian motions, we refer our readers to  \cite{ABBS, ABK}.

To maximize the possibility of $M_t\le \sqrt{2}\delta t$ for $\delta<1$, a good strategy
 is to make the first branching time $\tau$ as large as possible.
It was shown in \cite{CHM} that, conditioned on $\{M_t\le \sqrt{2}\delta t\}$, $\tau \approx \frac{1-\delta}{\sqrt{2}} t\pm O(1) \sqrt{t}$ when $\delta \in (1-\sqrt{2}, 1)$;
$\tau \approx t-O(1) \sqrt{t}$ when $\delta=1-\sqrt{2}$ and $\tau \approx t-O(1) $ when $\delta<1-\sqrt{2}$.
The asymptotic behaviors of  $P(M_t\le  \sqrt{2}\delta t)$ are different in these 3 different cases.

The intuition above also works for super-Brownian motion, but we need to use the first branching time of  the skeleton process, which is a branching Brownian motion.
Put
	 $$
q:=\psi'(\lambda^*)>0,
\qquad \rho:=\sqrt{1+\frac{\psi'(\lambda^*)}{\alpha}}=\sqrt{1+\frac{q}{\alpha}}.
$$
 We also use $\tau$ to denote the first branching time of
 the skeleton process of super-Brownian motion.
 We will prove that, conditioned on
 $\{M_t\le \sqrt{2\alpha}\delta t, \mathcal{S}\}$,
 as $t\to\infty$,
$ \tau\in [\frac{1-\delta}{\rho} t- (\log t)\sqrt{t}, \frac{1-\delta}{\rho} t+ (\log t)\sqrt{t}]$ when $ \delta\in (1-\rho,1)$;
$\tau \in t-\sqrt{t}\left[t^{-1/4}, \log t\right]$ when $ \delta=1-\rho$ and
$\tau\in [t-O(1),t] $ when $\delta<1-\rho$.
The asymptotic behavior of
$\P(M_t\le  \sqrt{2\alpha}\delta t|\mathcal{S})$
exhibits a phase transition at $\delta=1-\rho$.

Now we state our main results.
	
	\begin{theorem}\label{thm-case1}
		Assume that (H1) and (H2) hold.
		If $\delta\in(1-\rho, 1)$, then for any $f\in\mathcal{H}$,
		\begin{align*}
			&\lim_{t\to\infty}
			 e^{2\alpha(\rho-1)(1-\delta)t}t^{-3(\rho-1)/2}
			\E\left( e^{-\int_\R f(y-\sqrt{2\alpha}\delta t)X_t(dy)}; M_t\le \sqrt{2\alpha}\delta t|\mathcal{S}\right)\\
			=&\frac{\lambda^*}{e^{\lambda^*}-1}\frac{a_\delta^{3(\rho-1)/2}}{\sqrt{2\alpha}\rho}
			\int_{-\infty}^\infty e^{-\sqrt{2\alpha}(\rho-1)z}A(w_f(z))\d z,
		\end{align*}
		where $a_\delta=1-\frac{1-\delta}{\rho}$
		and
		$$
		A(\lambda)=\frac{1}{\lambda^* }\psi(\lambda)+\psi'(\lambda^*)
		\left(1-\frac{\lambda}{\lambda^*}\right)
		\ge0,\quad \lambda\geq 0.
		$$
	\end{theorem}

	\begin{theorem}\label{them:case2}
			Assume that (H1) and (H2) hold. Then
			for any $f \in\mathcal{H}$,
		\begin{align*}
			&\lim_{t\to\infty}t^{-3(\rho-1)/4}e^{(q+\alpha(\rho-1)^2)t}\E\left( e^{-\int_\R f(y-\sqrt{2\alpha}(1-\rho) t)X_t(dy)}; M_t\le \sqrt{2\alpha}(1-\rho) t|\mathcal{S}\right)\\
			=&\frac{\lambda^*}{e^{\lambda^*}-1}\frac{1}{\sqrt{2\pi }}
			\int^{\infty}_{ 0}s^{3(\rho-1)/2}e^{-\alpha\rho^2 s^2}\d s\int^\infty_{-\infty} e^{-\sqrt{2\alpha}(\rho-1)z}A(w_f(z))\d z.
		\end{align*}
	\end{theorem}
	
	\begin{theorem}\label{them:case3}
			Assume that (H1) and (H2) hold.
		If $\delta<1-\rho$, then for any $f \in\mathcal{B}_b^+(\R)$,
		\begin{align*}
  &\lim_{t\to\infty}\sqrt{t}e^{(q+\alpha\delta^2)t}\E\left( e^{-\int_\R f(y-\sqrt{2\alpha}\delta t)X_t(dy)}; M_t\le\sqrt{2\alpha}\delta t|\mathcal{S}\right)\\
=&\frac{\lambda^*}{e^{\lambda^*}-1}\left[\frac{1}{2\sqrt{\pi\alpha}|\delta|}
+\frac{1}{\sqrt{2\pi}}\int_0^\infty e^{(q-\alpha\delta^2)s}\,\d s\int_{\R}e^{\sqrt{2\alpha}\delta z}{G}_f(s,z)\,\d z\right],	
	\end{align*}
		where
\begin{equation}\label{e:defofG}
		G_f(t,x):=\frac{1}{\lambda^* }\Big[\psi(u_f(t,x))-\psi(\lambda^*+u^*_f(t,x))\Big]+qv_f(t,x),
\end{equation}
with $v_f,u^*_f$ being defined in
\eqref{def:v} and \eqref{def:ustar} below.
	\end{theorem}

The reason that we assume $f\in\mathcal{H}$ in Theorems \ref{thm-case1} and \ref{them:case2} is that \eqref{M-w} plays an important role in the proofs of Lemmas \ref{lem-case1} and \ref{lem-case2}. Lemma \ref{lem-case1}  is used in the proof of
Theorem  \ref{thm-case1} and Lemma \ref{lem-case2} is used in the proof of
Theorem  \ref{them:case2}.

Let ${\cal C}_c(\R)({\cal C}_c^+(\R))$ be the space of
all the (nonnegative) continuous functions with compact support.
Let $\cM_{R}(\R)$ be the space of
all the Radon measures on $\R$ equipped with the vague topology, see \cite[p.111]{Kallenberg}.
Recall that
for random measures $\mu_t,\mu\in \cM_{R}(\R)$,
$\mu_t$ converges
in distribution to $\mu$
is equivalent to $\langle f,\mu_t\rangle$ converges
in distribution to $\langle f,\mu\rangle $ for any $f\in{\cal C}_c^+(\R)$.
See \cite[p.119]{Kallenberg} for more details.

As a consequence of Theorems \ref{thm-case1}-\ref{them:case3}, we have the following corollary.

\begin{corollary}
	Assume that (H1) and (H2) hold.
	Conditioned on $\{M_t\le \sqrt{2\alpha}\delta t,\mathcal{S}\}$,
	$X_t-\sqrt{2\alpha}\delta t$ converges in distribution
 to a random measure $\Xi_\delta$. Moreover,
for any $f\in C_c^+(\R)$,
	if $\delta \in[1-\rho,1)$,
	\begin{align}\label{lap-xi}
		\E\left( e^{-\int_\R f(y)\Xi_\delta(dy)}\right)=\frac{	\int_{-\infty}^\infty e^{-\sqrt{2\alpha}(\rho-1)z}A(w_f(z))\d z}{	\int_{-\infty}^\infty e^{-\sqrt{2\alpha}(\rho-1)z}A(w(z))\d z};
	\end{align}
	and if $\delta<1-\rho$,
	$$\E\left( e^{-\int_\R f(y)\Xi_\delta(dy)}\right)=\frac{\frac{1}{\sqrt{2\alpha}|\delta|}
		+\int_0^\infty e^{(q-\alpha\delta^2)s}\,\d s\int_{\R}e^{\sqrt{2\alpha}\delta z}{G}_f(s,z)\,\d z}{\frac{1}{\sqrt{2\alpha}|\delta|}
		+\int_0^\infty e^{(q-\alpha\delta^2)s}\,\d s\int_{\R}e^{\sqrt{2\alpha}\delta z}{G}(s,z)\,\d z},$$
where $G_f$ is defined in \eqref{e:defofG}
and $G(t,x):=G_0(t,x)$.
\end{corollary}

{\bf Proof:}
First consider the case of  $\delta \in[1-\rho,1)$. For any $f\in\mathcal{H}$ and $\theta>0$, by Theorems \ref{thm-case1}-\ref{them:case2},
$$\lim_{t\to\infty}\E\left( e^{-\theta\int_\R f(y-\sqrt{2\alpha}\delta t)X_t(dy)}|M_t\le \sqrt{2\alpha}\delta t,\mathcal{S}\right)=\frac{	\int_{-\infty}^\infty e^{-\sqrt{2\alpha}(\rho-1)z}A(w_{\theta f}(z))\d z}{	\int_{-\infty}^\infty e^{-\sqrt{2\alpha}(\rho-1)z}A(w(z))\d z}.$$
It has been proved in \cite[Lemma 3.3]{SRZ} that $\lim_{\theta\to0}\tilde{C}(\theta f)=\tilde{C}(0)$, which implies that
$w_{\theta f}(x)\to w(x).$
Note that $A(\lambda)$ is decreasing on $(0,\lambda^*)$ and $0\le w_{\theta f}(z)\le \lambda^*$. Thus using
the monotone convergence theorem we get that
$$\lim_{\theta\to0}\frac{	\int_{-\infty}^\infty e^{-\sqrt{2\alpha}(\rho-1)z}A(w_{\theta f}(z))\d z}{	\int_{-\infty}^\infty e^{-\sqrt{2\alpha}(\rho-1)z}A(w(z))\d z}=1.$$
Thus, conditioned on $\{M_t\le \sqrt{2\alpha}\delta t,\mathcal{S}\}$,
$\int_\R f(y-\sqrt{2\alpha}\delta t)X_t(dy)$ converges in distribution for any $f\in \mathcal{C}_c^+(\R)$,
 which implies that $X_t-\sqrt{2\alpha}\delta t$
 converges in distribution
to a random measure $\Xi_\delta$ with Laplace transform given by \eqref{lap-xi}.

Similarly, using Theorem \ref{them:case3}, we can get the result for $\delta<1-\rho$.
\hfill$\Box$

\medskip

Throughout this paper we use $C$ to denote a positive constant whose value may change from one appearance to another.
For any two positive functions $f$ and $g$ on $[0,\infty)$,   $f\sim g$ as $s\to \infty$
 means that $\lim_{s\to\infty} \frac{f(s)}{g(s)}=1.$

	\section{Preliminaries}
	\subsection{Skeleton decomposition}\label{skeleton}
	
	Denote by $\P^*_{\mu}$ the law of $X$ with initial configuration $\mu$ conditioned on extinction.
	It is well known that
	$(X,\P^*)$ is a super-Brownian motion with branching mechanism $\psi^*(\lambda)=\psi(\lambda+\lambda^*)$. Note that $(\psi^*)'(0+)=\psi'(\lambda^*)=q>0$. So $(X,\P^*)$ is subcritical.

	Let $\mathbb{D}([0,\infty), \mathcal{M}_F({\R}))$ be the space of all the right continuous  functions $w:[0,\infty)\to \mathcal{M}_F({\R})$,  and $\mathbb{D}_0^+$ be the space of right continuous functions from $(0,\infty)$ to $ \mathcal{M}_F({\R})$ having zero as a trap. It has been proved in \cite{E.B2} that there is a family of measures $\{\N_x,x\in\R\}$ on $\mathbb{D}_0^+$ associated with the probability
	measures $\{\P_{\delta_x}^*:x\in\R\}$ such that
	\begin{equation}\label{N-measure}
		\int_{\mathbb{D}_0^+} \left(1-e^{-\langle f,w_t\rangle}\right)\N^*_{x}(dw)=-\log \P_{\delta_x}^*\left(e^{-\langle f,X_t\rangle}\right),
	\end{equation}
	for all $f\in\mathcal{B}_b^+(\R)$ and $t>0$.  The branching property of $X$ implies that, under $\P_{\delta_x}^*$, $X_t$ is an infinitely divisible measure, so \eqref{N-measure} is a Levy-Khinchine formula in which $\N^*_x$ plays the role of L\"evy measure.  By the spatial homogeneity of Brownian motion, one can check that
	$$\P_{\delta_x}^*\left(e^{-\langle f,X_t\rangle}\right)=\P_{\delta_0}^*\left(e^{-\int  f(x+y) X_t(\d y)}\right),\quad \N^*_{x}\left(1-e^{-\langle f,w_t\rangle}\right)=\N^*_{0}\left(1-e^{-\int  f(x+y) w_t(\d y)}\right).$$

	It was shown in \cite{BKS} that the skeleton of the super Brownian motion $X_t$  is a branching Brownian motion $Z_t$ with branching rate $q=\psi'(\lambda^*)$ and an offspring distribution $\{p_n:n\ge 2\}$ such that its generating function $\varphi$
	satisfies  $$q(\varphi(s)-s)=\frac{1}{\lambda^*}\psi(\lambda^*(1-s)).$$  We label the particles in $Z$ using the classical Ulam-Harris notation.  Let $\mathcal{T}$ be the set of all the particles. We write $\varnothing$ for the root.  For each particle $u\in\mathcal{T}$, we write $b_u$ and $\sigma_u$ for its birth and death time respectively, $N_u$ for the number of offspring of $u$, and $\{z_u(r):r\in[b_u,\sigma_u]$ for its spatial trajectory.  $v\preccurlyeq u$ means that $v$ is an ancestor of $u$.  Now we introduce the three kinds of immigrations along the skeleton $Z$ as follows.
	\begin{enumerate}
		\item {\bf Continuous immigration:} The process $I^{\N^*}$ is
		defined by
		$$
		I^{\N^*}_t:=\sum_{u\in{\cal T}}\sum_{(r_j,w_j)\in \mathcal{D}_{1,u}}{\bf 1}_{r_j<t} w_j(t-r_j),
		$$
		where, given $Z$, independently for each $u\in{\cal T}$, $\mathcal{D}_{1,u}:=\{(r_j,w_j):j\ge 1\}$ are the atoms of a
		Poisson point process on $(b_u, \sigma_u]\times
		\mathbb{D}_0^+$
		with rate
		$2\beta{\rm d} r\times{\rm d}\N^*_{z_u(r)}$.
		
		\item {\bf Discontinuous immigration:} The processes $I^{\P^*}$ is
		defined by
		$$
		I^{\P^*}_t:=\sum_{u\in{\cal T}}\sum_{(r_j,w_j)\in \mathcal{D}_{2,u}} {\bf 1}_{r_j<t}  w_j(t-r_j),
		$$
		where, given $Z$, independently for each $u\in{\cal T}$,
		$\mathcal{D}_{2,u}:=\{(r_j,w_j):j\ge 1\}$ are the atoms of a
		Poisson point process on $(b_u, \sigma_u]\times \mathbb{D}([0,\infty), \mathcal{M}_F({\R}))$ with rate
		${\rm d}r\times\int_{y\in(0,\infty)}ye^{-\lambda^* y}n({\rm }dy){\rm d}
		\P^*_{y\delta_{z_u(r)}}$.
		\item {\bf Branching point biased immigration:} The process $I^{\eta}$ is
		defined by
		\begin{equation*}
			I^{\eta}_t:=\sum_{u\in{\cal T}} \mathbf{1}_{\sigma_u\leq t}X^{(3,u)}_{t-\sigma_u}\ ,
		\end{equation*}
		where, given $Z$, independently for each $u\in{\cal T}$,
		$X^{(3, u)}_{\cdot}$ is an independent copy of the canonical process $X$ issued
		at time ${\sigma_u}$ with law $\mathbb{P}^*_{Y_u\delta_{z_u({\sigma_u})}}$ where,
		given $u$ has $n (\geq 2)$ offspring, $Y_u$ is an independent random variable with distribution
		$\eta_n({\rm d}y)$, where
		\begin{equation*}
			\eta_n(\d y)=\frac{1}{p_n\lambda^* q}\left\{\beta(\lambda^*)^2
			\delta_0(\d y)\mathbf{1}_{\{n=2\}}+(\lambda^*)^n\frac{y^n}{n!}e^{-\lambda^*y}n(\d y)\right\}.
		\end{equation*}
	\end{enumerate}
	Now we define another $\mathcal{M}_F(\mathbb{R})$-valued process $I=\{I_t : t\geq 0\}$ by
	\begin{equation}\label{I=sum}
		I:=I^{\N^*}+I^{\P^*}+I^{\eta}\ ,
	\end{equation}
	where $I^{\N^*}=\{I^{\N^*}_t: t\geq0\}$, $I^{\P^*}=\{I^{\P^*}_t: t\geq0\}$ and
	$I^{\eta}=\{I^{\eta}_t: t\geq0\}$,
	conditioned on $Z$, are independent of each other.
	For any integer-valued measure $\nu$, we denote by $\bP_{\nu}$ the law
	of $(Z, I)$ when the initial configuration of $Z$ is $\nu$.
	We write $\bP$ for  $\bP_{\delta_0}$.

	For any $\mu \in \mathcal{M}_F(\mathbb{R})$, let $Z$ be a branching Brownian motion with
	$Z_0$ being a Poisson random measure with intensity
	measure $\lambda^*\mu$ and $I$  is  the immigration process along $Z$. Let
	$\widetilde{X}$ be an independent copy of $X$ under $\P^*_\mu$,
	also independent of $I$.
	Then we define a measure-valued process $\Lambda=\{\Lambda_t: t\geq0\}$ by
	\begin{equation}\label{12}
		\Lambda=\widetilde{X}+I.
	\end{equation}
We denote  the law of $\Lambda$ by ${\bf Q}_\mu$.
	In particular, under $\bQ_{\delta_0}$, $Z_0=N\delta_0$, where $N$ is a Poisson random variable with parameter $\lambda^*$.
	We write $\bQ$ for $\bQ_{\delta_0}$.
In the rest of the paper,
we use ${\bf E}$,  $\mathbb{E}^*$ and $\mE_{\bQ}$ to denote the expectations with respect to ${\bf P}$, $\mathbb{P}^*$ and ${\bf Q}$, respectively.
	The following result is proved in \cite{BKS}.

	\begin{proposition}\label{p:skeleton}
		For any $\mu \in\mathcal{M}_F(\mathbb{R}^d)$, the process $(\Lambda, \bQ_\mu)$ is Markovian and has the same law as $(X, \P_\mu)$.
	\end{proposition}

	Recall that $M_t$ is the supremum of the support of $X_t$.
	Denote the supremum of
	$\Lambda_t, I_t, Z_t$, and $\tilde{X}_t$ by $M_t^\Lambda, M_t^I, M_t^Z$, and $ M_t^{\tilde{X}}$,  respectively.
By \eqref{V}, for any $f\in \mathcal{B}^+(\R)$,
$$
V_f(t,x)
=-\log \E_{\delta_x}\left(e^{-\int_\R f(y) X_t(dy)}\right), \quad x\in\R.$$ By the space homogeneity of $X$, we have
\begin{equation}\label{V(-x)}
V_f(t,-x)=-\log \E\left(e^{-\int_\R f(y-x) X_t(dy)}\right),\quad x\in\R.
\end{equation}
Setting $f_\theta:=f+\theta{\bf 1}_{(0,\infty)}$, we get
\begin{equation}\label{rel-V-u}
	u_f(t,x)=\lim_{\theta\to\infty} V_{f_\theta}(t,-x),\quad x\in\R.
\end{equation}
	For any $f\in \mathcal{B}^+(\R)$, put
\begin{align}
		v_f(t,x):= &\bE\left(e^{-\int_\R f(y-x) X_t(dy)}; M_t^{I}\le x\right),\label{def:v}\\
		u_f^*(t,x):=&-
		\log \E^*\left(e^{-\int_\R f(y-x) X_t(dy)}; M_t\le x\right).
		\label{def:ustar}
	\end{align}

For $f\equiv 0$, we write
$v(t,x)$ and  $u^*(t,x)$ for $v_f(t,x)$ and  $u^*_f(t,x)$, respectively.
The relation among $u_f, u^*_f$ and $v_f$ is given by the following lemma.
	
	\begin{lemma}\label{fact5}
For any $f\in \mathcal{B}^+(\R)$, $t\ge 0$ and $x\in\R$,
		$$u_f(t,x)=u^*_f(t,x)+\lambda^*(1-v_f(t,x)).$$
	\end{lemma}
	{\bf Proof:}
Recall that under $\bQ$, $Z_0=N\delta_0$, where $N$ is Poisson distributed with parameter $\lambda^*$.
By the definition of $\Lambda$, we get that,
for any $t\geq 0, x\in\R$,
	\begin{align*}
		e^{-u_f(t,x)}&=\E \left( e^{-\int_\R f(y-x) X_t(dy)}; M_t\le x\right)=\mE_\bQ \left( e^{-\int_\R f(y-x) \Lambda_t(dy)}; M_t^\Lambda\le x\right)\\
		&=\mE_\bQ \left( e^{-\int_\R f(y-x) \tilde{X}_t(dy)}; M_t^{\tilde{X}}\le x\right)\mE_\bQ \left( e^{-\int_\R f(y-x) I_t(dy)}; M_t^I\le x\right)\\
		&=\E^*\left(e^{-\int_\R f(y-x) X_t(dy)}; M_t\le x\right)
\mE_{\bQ}\left(\left[\bE(e^{-\int_\R f(y-x) I_t(dy)}; M_t^I\le x\right]^N\right)\\
&=e^{-u^*_f(t,x)}e^{\lambda^*(v_f(t,x)-1)}.
	\end{align*}
	Thus $u_f(t,x)=u^*_f(t,x)+\lambda^*(1-v_f(t,x)).$
	\hfill$\Box$
	\medskip

Now we give some basic relations among $M_t^Z,  M_t^\Lambda$ , $M_t^I$ and $M_t$.

	\begin{lemma}\label{fact1}
		Under $\bQ$, given $\Lambda_t$,
		$Z_t$ is a Poisson random measure with intensity $\lambda^*\Lambda_t$, which implies that $M_t^Z\le M_t^\Lambda$, $\bQ$-a.s.
	\end{lemma}
	{\bf Proof:}
	We refer the readers to  the display above \cite[(3.14)]{BKS} for a proof.
	\hfill$\Box$

	\begin{lemma}\label{fact2}
		Under $\bP$, $M_t^Z\le M_t^I$, a.s.
	\end{lemma}
	{\bf Proof:} First we claim that $\bQ(M_t^Z\le M_t^I)=1$. In fact, for any $x$, by
	Lemma \ref{fact1}, we have
		\begin{align*}
		0=\bQ(M_t^Z>x\ge  M_t^\Lambda)&=\bQ(M_t^Z>x, M_t^I\le x,M_t^{\tilde{X}}\le x)
		\\&
		=\bQ(M_t^Z>x, M_t^I\le x)
		\bQ(M_t^{\tilde{X}}\le x).
	\end{align*}
	Using the fact that $\bQ(M_t^{\tilde{X}}\le x)>0$, we get $\bQ(M_t^Z>x, M_t^I\le x)=0$.
	Since $x$ is arbitrary, the claim is true.
	
	Recall that under $\bQ$, $Z_0=N\delta_0$, where $N$ is Poisson distributed with parameter $\lambda^*$. Thus
	$$0=\bQ(M_t^Z> M_t^I)\ge \bQ[M_t^Z>M_t^I|N=1]\bQ(N=1)=\bP(M_t^Z> M_t^I)e^{-\lambda^*},$$
	which implies that $\bP(M_t^Z>M_t^I)=0$.
	\hfill$\Box$

	\medskip

	The following lemma implies that,
	to prove our main results,  we only need to study the limit behavior of  $v_f(t,\sqrt{2\alpha}\delta t).$
	
	\begin{lemma}\label{fact4}
For any $f\in \mathcal{B}^+(\R)$ and $\delta<1$,
		\begin{equation}\label{equalent}
           \lim_{t\to\infty} \frac{\E\left(e^{-\int_\R f(y-\sqrt{2\alpha}\delta t) X_t(dy)}; M_t\le \sqrt{2\alpha}\delta t|\mathcal{S}\right)}{v_f(t,\sqrt{2\alpha}\delta t)}=\frac{\lambda^*}{e^{\lambda^*}-1}.
		\end{equation}
	\end{lemma}
	{\bf Proof:}
We also use $\mathcal{S}$ to denote the survival of $\Lambda$. It is clear that, under $\bQ$, $\mathcal{S}\subset \{N\ge 1\}$  and $\bQ(\mathcal{S})=\bQ(N\ge1)=1-e^{-\lambda^*}$. It follows that
	$\mathcal{S}=\{N\ge 1\}$, $\bQ$-a.s.
	Then, by Proposition \ref{p:skeleton},
\begin{align}\label{e:1}
	&\E\left(e^{-\int_\R f(y-x) X_t(dy)}; M_t\le x|\mathcal{S}\right)=\mE_\bQ \left( e^{-\int_\R f(y-x) \Lambda_t(dy)}; M_t^\Lambda\le x|N\ge 1\right)\nonumber\\
	=&\mE_\bQ \left( e^{-\int_\R f(y-x) \tilde{X}_t(dy)}; M_t^{\tilde{X}}\le x\right)\mE_\bQ \left( e^{-\int_\R f(y-x) I_t(dy)}; M_t^I\le x|N\ge 1\right)\nonumber\\
	=&e^{-u^*_f(t,x)}\mE_\bQ(v_f(t,x)^N|N\ge 1)
	=e^{-u^*_f(t,x)}\frac{e^{\lambda^*v_f(t,x)}-1}{e^{\lambda^*}-1 }.
\end{align}
	Since $(X,\P^*)$ is subcritical,
	we have, for any $\delta$, $$e^{-u^*_f(t,\sqrt{2\alpha}\delta t)}\ge \P^*(\|X_t\|=0)\to 1,\quad t\to\infty,$$ which implies that $e^{-u^*_f(t,\sqrt{2\alpha}\delta t)}\to1$, as $t\to\infty$.
	By \eqref{limit-as}, we have for any $\delta<1$,
	$$
	\E\left(e^{-\int_\R f(y-\sqrt{2\alpha}\delta t) X_t(dy)}; M_t\le \sqrt{2\alpha}\delta t|\mathcal{S}\right)\le \P(M_t\le \sqrt{2\alpha}
	\delta t|\mathcal{S})\to0.
	$$
	Thus by \eqref{e:1},
	$v_f(t,\sqrt{2\alpha}\delta t)\to0 $ for any $\delta<1$. The desired result follows immediately.
	\hfill $\Box$
	
	\medskip

To study the behavior of $v_f(t,\sqrt{2\alpha}\delta t)$ as $t\to\infty$, the following decomposition of $v_f$ plays a fundamental role.

	\begin{proposition}\label{exp:v}
For any $f\in \mathcal{B}^+(\R)$,  $t>0$ and $x\in\R$,
\begin{equation}\label{decom-v}
v_f(t,x)=U_{1,f}(t,x)+U_{2,f}(t,x),
\end{equation}
where
		\begin{align}\label{eq-v2}
			U_{1,f}(t,x)=&\mE
			\Big[e^{-\int_0^t
				\psi'(\lambda^*+u^*_f(t-r,x-B_r))\,dr},B_t\le x\Big],\\
			U_{2,f}(t,x)=&
			\mE\int_0^t e^{-\int_0^{s} \psi'(\lambda^*+u^*_f(t-r,x-B_{r}))\,dr}\hat{G}_f(t-s,x-B_{s})\,\d s,\label{eq-v2'}
		\end{align}
with $\hat{G}_f(t,x)$ being defined by
		\begin{align*}
			\hat{G}_f(t,x)&=\frac{1}{\lambda^*}\left[\beta(\lambda^*)^2 v_f(t,x)^2+\int_0^\infty \left(e^{\lambda^* v_f(t,x)y}-1-\lambda^* v_f(t,x)y\right)e^{-(\lambda^*+u^*_f(t,x))y}\,n(\d y)\right]\\
			&=\frac{1}{\lambda^*}\Big[\psi(u_f(t,x))-\psi(\lambda^*+u^*_f(t,x))+\psi'(\lambda^*+u^*_f(t,x))\lambda^*v_f(t,x)\Big].
		\end{align*}
	\end{proposition}
	{\bf Proof:}
Let $\tau$ be the first splitting time of $Z$, that is $\tau=\sigma_\varnothing$.
	By considering the cases $\tau>t$ and $\tau\le t$ separately, we get
	\begin{align}\label{v=U1+U2}
		&v_f(t,x)=\bE\left( e^{-\int_\R f(y-x)I_t(dy)}; M^I_t\le x\right)\nonumber\\
		&=\displaystyle \bE\left( e^{-\int_\R f(y-x)I_t(dy)}; M^I_t\le x,\tau>t\right)
		+\bE\left( e^{-\int_\R f(y-x)I_t(dy)}; M^I_t\le x,\tau\le t\right)\nonumber\\
&=:\displaystyle U_{1,f}(t,x)+U_{2,f}(t,x).\end{align}

By Lemma \ref{fact2},
	$U_{1,f}(t,x)=\bE\left( e^{-\int_\R f(y-x)I_t(dy)};  M^I_t\le x, M_t^Z\le x,\tau>t\right)$.
	By the decomposition of $I$ in
	\eqref{I=sum},
	on the event $\{\tau>t\}$, we have that $I_t=I_t^{\N^*}+I_t^{\P^*}$.  Thus using \cite[Lemma 3]{BKS},
	we have on the event $\{\tau>t\}$, for any $x\in\R$,
	\begin{align*}
		&\bE\left( e^{-\int_\R f(y-x)I_t(dy)}; M^I_t\le x|\mathcal{F}_t^Z\right)
   =\lim_{\theta\to\infty}\bE\left(e^{-\int_\R [f(y-x)+\theta {\bf 1}_{(0,\infty)}(y-x)]I_t(dy)}|\mathcal{F}_t^Z\right)\\
		=&\exp\left\{-\int_0^t \langle\phi(u^*_f(t-s, x-\cdot)),Z_s\rangle\,\d s\right\},
	\end{align*}
	where
$\{\mathcal{F}_t^Z, t\geq 0\}$
is the natural filtration of $Z$ and
	\begin{equation}\label{phi}
		\phi(\lambda):=\psi'(\lambda+\lambda^*)-\psi'(\lambda^*)=2\beta\lambda+\int_0^\infty (1-e^{-\lambda x})xe^{-\lambda^*x}n(dx).
	\end{equation}
	Note that, on the event $\{\tau>t\}$, $Z_s=\delta_{z_{\varnothing}(s)}$ and $\{z_{\varnothing}(s),s\le t\}\overset{d}{=}\{B_s,s\le t\}$.  Thus
	\begin{align}\label{U1-expansion}
		U_{1,f}(t,x)&=e^{-qt}\mE\Big[\exp\Big\{-\int_0^t \phi(u^*_f(t-r,x-B_r))\,dr\Big\}; B_t\le x\Big]\nonumber\\
		&=\mE\Big[\exp\Big\{-\int_0^t \psi'(\lambda^*+u^*_f(t-r,x-B_r))\,dr\Big\};B_t\le x\Big].
	\end{align}

	On the event $\{\tau\le t\}$, the immigration process $I$ has the following
	expression:
	\begin{align}\label{dec:I}
		I_t&
		=\sum_{(r_j,w_j)\in\mathcal{D}_{1,\varnothing}} w_j(t-r_j)+\sum_{(r_j,w_j)\in\mathcal{D}_{2,\varnothing}} w_j(t-r_j)+
		X^{(3,\varnothing)}_{t-\tau}
		+\sum_{i=1}^{N_\varnothing}I^i_{t-\tau}\nonumber\\
		&=:
		\mathcal{J}_{1,t}+\mathcal{J}_{2,t}+\mathcal{J}_{3,t}+\mathcal{J}_{4,t},
	\end{align}
	where, given $Z_{\tau}$, $I^i,i=1,\cdots, N_\varnothing$, are i.i.d copies of $I$ under $\bP_{z_\varnothing(\tau)}$.
	Since, given $\mathcal{F}_t^Z$, $\mathcal{J}_{i,t},i=1,2,3,4$, are independent, so
	\begin{align}\label{2.5}
U_{2,f}(t,x)&=\bE\left[\bE\left( e^{-\int_\R f(y-x)I_t(dy)}; M^I_t\le x | \mathcal{F}_t^Z\right);\tau\le t\right]\nonumber\\
&=\bE\left[H_{1,t}H_{2,t}H_{3,t}H_{4,t};\tau\le t\right],
	\end{align}
where
$$H_{i,t}=\bE( e^{-\int_\R f(y-x)\mathcal{J}_{i,t}(dy)}; \mathcal{J}_{i,t}(x,\infty)=0 | \mathcal{F}_t^Z),\quad i=1,2,3,4.$$
Put $f_\theta=f+\theta{\bf 1}_{(0,\infty)}.$
By the bounded convergence theorem,
we have
\begin{equation}\label{H}
	H_{i,t}=\lim_{\theta\to\infty} \bE\left( e^{-\int_\R f_\theta(y-x)\mathcal{J}_{i,t}(dy)}| \mathcal{F}_t^Z\right).
\end{equation}
	By the definition of $\mathcal{D}_{1,\varnothing}$ and \eqref{H}, we have that, on the event $\{\tau\le t\}$,
	\begin{align}
		H_{1,t}=\lim_{\theta\to\infty} \exp\left\{-2\beta\int_0^\tau\int_{\mathbb{D}_0^+}\Big(1-e^{-\int_\R f_\theta(y-x) w_{t-r}(dy)}\Big)\N^*_{z_{\varnothing}(r)}(dw)\,dr\right\}.
	\end{align}
	Using \eqref{N-measure}, we getthat
	\begin{align*}
		&\lim_{\theta\to\infty}\int_{\mathbb{D}_0^+}\left(1-e^{-\int_\R f_\theta(y-x) w_{t-r}(dy)}\right)\N^*_{z}(dw)
		=\lim_{\theta\to\infty}-\log \E^*_{\delta_z}\left[e^{-\int_\R f_\theta(y-x) X_{t-r}(dy)}\right]\\
		&=-\log \E^*_{\delta_z}\left[e^{-\int_\R f(y-x) X_{t-r}(dy)};M_{t-r}\le x\right]=u^*_f(t-r,x-z).
	\end{align*}
	Thus we have that
	\begin{align}\label{J1}
		H_{1,t}=\exp\left\{-2\beta\int_0^\tau u^*_f(t-r,x-z_{\varnothing}(r))\,dr\right\}.
	\end{align}
	For $H_{2,t}$,  on the event $\{\tau\le t\}$, we have that
	\begin{align}
		H_{2,t}=\lim_{\theta\to\infty}\exp\left\{-\int_0^\tau \,dr \int_0^\infty ye^{-\lambda^* y}n(\d y)\E^*_{y\delta_{z_\varnothing}(r)}\Big(1-e^{-\int_\R f_\theta(y-x) X_{t-r}(dy)}\Big)\,dr\right\}.
	\end{align}
	It follows from the branching property of $X$ that
	\begin{align*}
	&\lim_{\theta\to\infty}\P^*_{y\delta_{z}}(e^{-\int_\R f_\theta(y-x) X_{t-r}(dy)})
		=\lim_{\theta\to\infty}\left[\P^*_{\delta_{z}}(e^{-\int_\R f_\theta(y-x) X_{t-r}(dy)})\right]^y=e^{-u^*_f(t-r,x-z)y},\end{align*}
	which implies that
	\begin{align}\label{J2}
		H_{2,t}=\exp\left\{-\int_0^\tau \int_0^\infty y[1-e^{-u^*_f(t-r,x-z_\varnothing(r))y}]e^{-\lambda^* y}n(\d y)\,dr\right\}.
	\end{align}
	By the definition of $X^{(3,\varnothing)}$, on the event $\{\tau\le t\}$, we have that
	\begin{align}\label{J3}
H_{3,t}
=&\lim_{\theta\to \infty}\bE\left(\P^*_{Y_\varnothing\delta_{y}}\left(e^{-\int_\R f_\theta(y-x)X_{t-s}(dy)}\right)|\mathcal{F}_t^Z\right)|_{s=\tau,y=z_\varnothing(\tau)}\nonumber\\
=&\bE\left(e^{-u^*_f(t-\tau,x-z_\varnothing(\tau))Y_\varnothing}|\mathcal{F}_t^Z\right)\nonumber\\
=&\frac{1}{p_{N_\varnothing}\lambda^* q}\left(\beta(\lambda^*)^2{\bf 1}_{N_\varnothing=2}+\int_0^\infty\frac{(\lambda^*y)^{N_\varnothing}}{N_\varnothing !}e^{-u^*_f(t-\tau,x-z_\varnothing(\tau))y}e^{-\lambda^* y}\,n(\d y)\right).
	\end{align}
	It follows from the branching property that
	on the event $\{\tau\le t\}$,
	\begin{align}\label{J4}
		H_{4,t}
		=\left[\bP_{\delta_{z_\varnothing(\tau)}}\Big(e^{-\int_\R f(y-x)X_{t-s}(dy)}; M_{t-s}^I\le x\Big)\right]^{N_{\varnothing}}_{s=\tau}=v_f(t-\tau,x-z_\varnothing(\tau))^{N_{\varnothing}}.
	\end{align}
	Note that
	\begin{align}\label{2.3}
		&\sum_{n=2}^\infty p_n\frac{1}{p_{n}\lambda^* q}\left(\beta(\lambda^*)^2{\bf 1}_{n=2}+\int_0^\infty\frac{(\lambda^*y)^{n}}{n!}e^{-u^*_f(t-\tau,x-z_\varnothing(\tau))y}e^{-\lambda^* y}\,n(\d y)\right)v_f(t-\tau,x-z_\varnothing(\tau))^{n}\nonumber\\
		&=\frac{1}{\lambda^* q}\Big[(\beta(\lambda^*)^2v_f(t-\tau,x-z_\varnothing(\tau))^2\nonumber\\
		&\qquad +\int_0^\infty\left(e^{\lambda^*v_f(t-\tau,x-z_\varnothing(\tau))y}-1-\lambda^*v_f(t-\tau,x-z_\varnothing(\tau))y\right)
e^{-(\lambda^*+u^*_f(t-\tau,x-z_\varnothing(\tau)))y}n(\d y)\Big]\nonumber\\
		&=q^{-1}\hat{G}_f(t-\tau,x-z_\varnothing(\tau)).
	\end{align}
	Recall the definition of $\phi$ in \eqref{phi}.  Combining \eqref{2.5}-\eqref{2.3}, we get that
	\begin{align*}
		U_{2,f}(t,x)&=q^{-1}\bP\left(\exp\left\{-\int_0^{\tau} \phi( u^*_f(t-r,x-z_\varnothing(r)))\,dr\right\}\hat{G}_f(t-\tau,x-z_\varnothing(\tau)),\tau\le t\right)\nonumber\\
		&=\mE\int_0^t\exp\left\{-\int_0^{s}\left(q+ \phi( u^*_f(t-r,x-B_{r}))\right)\,dr\right\}\hat{G}_f(t-s,x-B_{s})\,\d s.
	\end{align*}
	Note that $q+\phi(\lambda)=\psi'(\lambda^*+\lambda)$.
         The proof is now compete.
	\hfill$\Box$

	\medskip
	
	Note that $\frac{e^x-1-x}{x^2}=\sum_{k=2}^\infty \frac{x^{k-2}}{k!}$ is increasing in  $x$ on $(0,\infty)$. So
	$e^{\lambda^*v_f(t,x)y}-1-\lambda^*v_f(t,x)y\le v_f(t,x)^2(e^{\lambda^*y}-1-\lambda^*y)$, which implies that
	\begin{align}\label{est-hatG}
		&\hat{G}_f(t,x)\le \frac{1}{\lambda^*}\left[(\beta(\lambda^*)^2+\int_0^\infty(e^{\lambda^*y}-1-\lambda^*y)e^{-\lambda^*y}n(\d y))\right]v_f(t,x)^2
		\nonumber\\
		&=\left(\psi'(\lambda^*)-\psi(\lambda^*)/\lambda^*\right)v_f(t,x)^2=qv_f(t,x)^2
		\le qv(t,x).
	\end{align}
	Here in the last inequality, we use the fact that $v_f(t,x)\le v(t,x)$.

	\subsection{Some useful estimates}
	In this subsection we give
	some useful estimates for $u^*_f(t,x)$ and $v_f(t,x)$. Recall that $q=\psi'(\lambda^*)$ and $\rho=\sqrt{1+q/\alpha}.$

	\begin{lemma}\label{lemma:u*}
		\begin{description}
			\item  [{\bf (1)}]
			For any $f\in\mathcal{B}^+(\R)$ and $t>0, x\in\R$,
			$$u^*_f(t,x)\le k(t):=-\log \P^*(X_t=0),$$ and
			$t\mapsto  e^{qt}k(t)$ is decreasing on $(0, \infty)$.
		\item [{\bf (2)}]
			If (H2) holds, then there exists a positive constant
			$c_2$ such that
			\begin{equation}\label{est:k}
				k(t)\le  \left[\frac{c_2}
				{e^{c_2\vartheta t}-1}\right]^{1/\vartheta}, \quad t>0,
			\end{equation}
			and for any $f\in\mathcal{B}_b^+(\R)$,
			there exists a positive constant $c_3$ such that
			\begin{align}\label{est:u*}
				u^*_f(t,x)\le c_3(1+x^{-2/\vartheta})e^{(a+\alpha)t},
				\quad t, x>0.
			\end{align}
		\end{description}
	\end{lemma}
	{\bf Proof:}
	Since $\E^*\left(e^{-\int_R f(y-x)X_t(dy)}; M_t\le x\right)\ge \P^*(X_t=0)$
	for any $t>0, x\in\R$,
	we have $u^*_f(t,x)\le k(t)$.
	By the branching property and Markov property, we get that
	$$\P^*(\|X_t\|=0)=\E^*\left(\P^*_{X_{t-s}}(\|X_{s}\|=0)\right)=\E^*\left(e^{-k(s)\|X_{t-s}\|}\right).$$
	Put $u^*_\theta(t):=-\log\E^*\left(e^{-\theta\|X\|_t}\right)$. Then
	$k(t)=u^*_{k(s)}(t-s).$
	Under $\P^*$, $\|X_t\|$ is a continuous state branching process with branching mechanism $\psi(\lambda^*+\lambda)$. Then according to \cite[Theorem 10.1]{Kyprianou}, we have
\begin{equation}\label{deriv-k}
	k'(t)=-\psi\left(\lambda^*+u^*_{k(s)}(t-s)\right)=-\psi(\lambda^*+k(t)).
\end{equation}
Since $\psi(\lambda^*)=0$ and $\psi'$ is increasing on $(0,\infty)$,
 $\psi(\lambda^*+\lambda)\ge \psi'(\lambda^*)\lambda=q\lambda.$ Thus $k'(t)\le -q k(t)$.
	Using this one can check that $(e^{qt}k(t))'\le 0$.
	The proof of (1) is complete.
	
	Assume that (H2) holds.
	Then  there exists $c_2>0$ such that
	$\psi(\lambda^*+\lambda)\ge c_2(\lambda+\lambda^{1+\vartheta})$.
	Now \eqref{est:k} follows immediately from
\eqref{deriv-k}.
	Since $u_f^*(t,x)\le u_f(t,x)$, it suffices to show that \eqref{est:u*} is true for $u_f(t,x).$
	By \cite[Lemma 2.3(2)]{SRZ}, we have that
	$$V_{f_1+f_2}(t,x)\le V_{f_1}(t,x)+V_{f_2}(t,x).$$
By \eqref{rel-V-u},
$$u_f(t,x)=\lim_{\theta\to\infty} V_{f_\theta}(t,-x)\le V_f(t,-x)+\lim_{\theta\to\infty}V_{\theta{\bf 1}_{(0,\infty)}}(t,-x)
	=V_f(t,-x)+u(t,x),$$
where $f_\theta=f+\theta{\bf 1}_{(0,\infty)}$.
     By \eqref{V(-x)} and Jensen's inequality,  we have that
	$$V_f(t,-x)=-\log \E\left(e^{-\int_\R f(y-x) X_t(dy)}\right)\le \E\left(\int_\R f(y-x) X_t(dy)\right)=e^{\alpha t}\mE(f(B_t-x))\le e^{\alpha t}\|f\|.$$
	By \cite[Lemma 4.2 and 4.3]{SRZ} (with $A$ being replaced by $x$, and $x$ there replaced by 0),
	we get that there
	exists a positive constant
	$C$ such that
	$$
	u(t,x)\le
	C(1+x^{-2/\vartheta})e^{at},\quad t,x>0                                       .
	$$
Combining the two displays above, we get that
	$$u_f(t,x)\le e^{\alpha t}\|f\|+	C(1+x^{-2/\vartheta})e^{at}\le (C+\|f\|)(1+x^{-2/\vartheta})e^{(a+\alpha)t}.$$
	Now \eqref{est:u*} follows immediately.
	\hfill$\Box$
	
	\medskip
	
\begin{lemma}\label{lemma2}
Assume that (H1) and (H2) hold.		For any $A>0$ and  $\epsilon>0$,
		$$\int_0^A \phi(k(s))s^\epsilon\,\d s<\infty.$$
	\end{lemma}
{\bf Proof:}
Note that, by \eqref{deriv-k},
$$k'(s)=-\psi(k(s)+\lambda^*), \quad k''(s)=-\psi'(k(s)+\lambda^*)k'(s). $$
Thus, using \eqref{phi}, we have
$$0\le \phi(k(s))=\psi'(k(s)+\lambda^*)-q=\frac{k''(s)}{-k'(s)}-q\le \frac{k''(s)}{-k'(s)}.$$
It follows that
\begin{align}\label{3.10.1}
\int_0^A \phi(k(s))s^{\epsilon}\,\d s&\le \int_0^A \frac{k''(s)}{-k'(s)}s^{\epsilon}\,\d s=\int_0^A  s^{\epsilon}\,\d (-\log (-k'(s)))\nonumber\\
&=-\log(-k'(A))A^{\epsilon}+\lim_{s\to0} s^{\epsilon}\log (-k'(s))+\epsilon\int_0^A \log (-k'(s))s^{\epsilon-1}\,\d s.
\end{align}
Note that for $\lambda>0$,  $\psi''(\lambda+\lambda^*)$ exists and is decreasing. By Taylor's expansion, since $\psi(\lambda^*)=0$,  we have that
$$\psi(\lambda+\lambda^*)\le\psi'(\lambda^*)\lambda+\psi''(\lambda^*)\lambda^2,\quad \lambda>0.$$
By \eqref{est:k}, we have that $k(s)\le Cs^{-1/\vartheta}$.
Thus we get that
\begin{align*}
-k'(s)=\psi(k(s)+\lambda^*)&\le \psi'(\lambda^*)k(s)+\psi''(\lambda^*)k(s)^2\\
&\le C(s^{-1/\vartheta}+s^{-2/\vartheta})\le C s^{-2/\vartheta},\quad s\in[0,A].
\end{align*}
Now the desired result follows immediately from \eqref{3.10.1}.
\hfill$\Box$

\medskip
	
	Now we give some upper estimates of $v(t,x)$.

\begin{lemma}\label{lemma:key}
		\begin{description}
			\item [{\bf(1)}]For any $t>0$,
			\begin{equation}\label{domi-B'}
				v(t,x)\le {\rm P}(B_t\le x), \qquad x\in \R,
			\end{equation}
           and
			\begin{equation}\label{domi-B}
v(t,x)\le {\rm P}(B_t\le x) \le
\frac{\sqrt{t}}{\sqrt{2\pi}|x|} e^{-\frac{x^2}{2t}}, \qquad  x<0.
		\end{equation}
		
           \item [{\bf (2)}]
			There exist $t_0>1$ and $c>0$ such that for any $t>t_0$,
\begin{align}\label{est:v2}
  v(t,\sqrt{2\alpha} \theta t-\sqrt{t})&\le {\bP}({M}^Z_t\le \sqrt{2\alpha}  \theta t-\sqrt{t})\nonumber\\
   & \le c t\left\{\begin{array}{ll}
					e^{-(q+\alpha  \theta^2)t},&\hbox{$ \theta<1-\rho$};\\
					e^{-2\alpha(\rho-1)(1- \theta)t},&\hbox{$1-\rho\le  \theta<1$}.
				\end{array}\right.
\end{align}
		\end{description}
	\end{lemma}
{\bf Proof:} (1)
	By Proposition \ref{exp:v},
 we have
\begin{align}\label{eq-v2''}
v(t,x)=&\mE
	\Big[e^{-\int_0^t
		\psi'(\lambda^*+u^*(t-r,x-B_r))\,dr},B_t\le x\Big]\nonumber\\
	&+\mE\int_0^t e^{-\int_0^{s} \psi'(\lambda^*+u^*(t-r,x-B_{r}))\,dr}\hat{G}(t-s,x-B_{s})\,\d s\nonumber\\
	=&\mE_x
	\Big[e^{-\int_0^t
		\psi'(\lambda^*+u^*(t-r,B_r))\,dr},B_t\ge 0\Big]\nonumber\\
	&+\mE_x\int_0^t e^{-\int_0^{s} \psi'(\lambda^*+u^*(t-r,B_{r}))\,dr}\hat{G}(t-s,B_{s})\,\d s,
\end{align}
where $\hat{G}$ is the $\hat{G}_f$ defined in
Proposition \ref{exp:v} with $f\equiv 0$.
Thus, by the Feynman-Kac formula, we have
	\begin{align*}
		v(t,x)&=\mP_x \Big[B_t\ge  0\Big]+\mE_x \int_0^t\left[\hat{G}(t-s,B_{s})-\psi'(\lambda^*+u^*(s,B_{t-s}))v(s,B_{t-s})\right]\d s\\
		&=\mP_x \Big[B_t\ge 0\Big]+\frac{1}{\lambda^*}\mE_x \int_0^t\left[\psi(u(s, B_{t-s}))-\psi(\lambda^*+u^*(s,B_{t-s}))\right]\d s.
	\end{align*}
Note that $u(s,z)\le \lambda^*+u^*(s,z)$,   $\psi$ is negative on $(0,\lambda^*)$ and increasing on $(\lambda^*,\infty)$. Thus $\psi(u(s, B_{t-s}))-\psi(\lambda^*+u^*(s,B_{t-s}))\le 0.$ Therefore we have that
$$v(t,x)\le \mP_x \Big[B_t\ge 0\Big]=\mP \Big[B_t\le x\Big],\qquad x\in\R.$$

For  $x<0$,
\begin{align}\label{est:bm}
	\mP \Big[B_t\le x\Big]&=\mP \Big[B_1\ge |x|t^{-1/2}\Big]=\frac{1}{2\pi}\int_{|x|t^{-1/2}}^\infty e^{-y^2/2}\d y\nonumber\\
	&\le \frac{1}{2\pi }\int_{|x|t^{-1/2}}^\infty \frac{y}{|x|t^{-1/2}}e^{-y^2/2}\d y\le \frac{\sqrt{t}}{\sqrt{2\pi}|x|} e^{-x^2/(2t)}.
\end{align}
Thus \eqref{domi-B} follows.

	(2) We claim that there exists $t_0>0$ such that for any $t>t_0$ and $x$,
	\begin{align}\label{upper}
		{\bP}({M}^Z_t\le z)\le (2qt+1)\sup_{0\le s\le t}e^{-qs}\mP\left(B_{1}\le ( z-\sqrt{2\alpha}(t-s)+\sqrt{t})/\sqrt{s}\right).
	\end{align}
    It is shown in \cite{DS} (see the discussion below \cite[Lemma 3]{DS}) that the claim is true when $p_2=1$ and $q=1$. Using similar arguments we see that it is also true for the general case. We omit the proof here.
	
	Put $a(t):=\sqrt{2\alpha}(1- \theta)t.$
By \eqref{upper},  for $t>t_0$,
	$${\bP}({M}^Z_t\le \sqrt{2\alpha}  \theta t-\sqrt{t})\le (2qt+1)\sup_{0\le s\le t}e^{-qs}\mP(B_{1}\le (\sqrt{2\alpha}s-a(t))/\sqrt{s}).$$
	Note that by \eqref{est:bm}, $P(B_1\le -y)\le  \frac{1}{\sqrt{2\pi}} y^{-1}e^{-y^2/2}$ for all $y>0$. Thus, if $\sqrt{2\alpha}s<a(t)$, we have
	\begin{align}\label{6.1}
		e^{-qs}\mP\left(B_{1}\le (\sqrt{2\alpha}s-a(t))/\sqrt{s}\right)
		&\le \frac{\sqrt{s}}{\sqrt{2\pi}}\frac{1}{{a(t)-\sqrt{2\alpha}s}}e^{-qs} e^{-(\sqrt{2\alpha}s-a(t))^2/2s}\nonumber\\
		&=\frac{\sqrt{s}}{\sqrt{2\pi}}\frac{1}{{a(t)-\sqrt{2\alpha}s}}e^{\sqrt{2\alpha}a(t)}e^{-\alpha\rho^2s-\frac{a(t)^2}{2s}} .
	\end{align}
	It is clear that
\begin{equation}\label{easy-domi-below}
\alpha\rho^2s+\frac{a(t)^2}{2s}\ge \sqrt{2\alpha}\rho a(t),
\end{equation}
and is decreasing on $(0,\frac{a(t)}{\sqrt{2\alpha}\rho})$.
We now  prove the desired result in four cases.

	(i) If $a(t)>\sqrt{2\alpha}\rho t$ (that is, $ \theta<1-\rho$), then $\sqrt{2\alpha}s<a(t)$ for $s\in[0,t]$ and thus by \eqref{6.1} we have that
	\begin{align*}
		&\sup_{0\le s\le t}e^{-qs}\mP\left(B_{1}\le (\sqrt{2\alpha}s-a(t))/\sqrt{s}\right)
		\le \frac{\sqrt{t}}{\sqrt{2\pi}}\frac{1}{{a(t)-\sqrt{2\alpha}t}}e^{\sqrt{2\alpha}a(t)}e^{-\alpha\rho^2t-\frac{a(t)^2}{2t}}\\
		&= \frac{\sqrt{t}}{\sqrt{2\pi}}\frac{1}{{a(t)-\sqrt{2\alpha}t}}e^{-(\alpha\rho^2+\alpha(1- \theta)^2-2\alpha(1- \theta))t}
		\le\frac{1}{\sqrt{2\pi}}\frac{1}{\sqrt{2\alpha}(\rho-1)\sqrt{t}}e^{-(q+\alpha  \theta^2)t}.
	\end{align*}

	(ii) If $\sqrt{2\alpha}\frac{\rho+1}{2} t \le  a(t)\le\sqrt{2\alpha}\rho t$ (that is, $1-\rho\le  \theta\le (1-\rho)/2$), then  $\sqrt{2\alpha}s<a(t)$ for $s\in[0,t]$, and thus by
\eqref{6.1} and \eqref{easy-domi-below} we have that
	\begin{align*}
		\sup_{0\le s\le t}e^{-qs}\mP\left(B_{1}\le (\sqrt{2\alpha}s-a(t))/\sqrt{s}\right)
		&\le\frac{1}{\sqrt{2\pi}}\frac{2}{\sqrt{2\alpha}(\rho-1)\sqrt{t}}e^{-\sqrt{2\alpha}(\rho-1)a(t)}\\
		&=\frac{1}{\sqrt{2\pi}}\frac{2}{\sqrt{2\alpha}(\rho-1)\sqrt{t}}e^{-2\alpha(\rho-1)(1- \theta)t}.
	\end{align*}

	(iii) If $1<a(t)<\sqrt{2\alpha}\frac{\rho+1}{2}t$ (that is, $(1-\rho)/2< \theta<1-\frac{1}{\sqrt{2\alpha}t}$), then
	\begin{align*}
		&\sup_{0\le s\le t}e^{-qs}\mP\left(B_{1}\le (\sqrt{2\alpha}s-a(t))/\sqrt{s}\right)\\
		\le & \sup_{0\le s\le \frac{2}{\sqrt{2\alpha}(\rho+1)}a(t)}e^{-qs}\mP\left(B_{1}\le (\sqrt{2\alpha}s-a(t))/\sqrt{s}\right)+e^{-q\frac{2}{\sqrt{2\alpha}(\rho+1)}a(t)}\\
		\le& \sup_{0\le s\le \frac{2}{\sqrt{2\alpha}(\rho+1)}a(t)}\frac{1}{\sqrt{2\pi}}\frac{\sqrt{s}}{a(t)-\sqrt{2\alpha}s}e^{-\sqrt{2\alpha}(\rho-1)a(t)}+e^{-\sqrt{2\alpha}(\rho-1)a(t)}\\
		\le  & \sqrt{\frac{1}{\sqrt{2\alpha}\pi(\rho+1)}}\frac{1}{\frac{\rho-1}{(\rho+1)}\sqrt{a(t)}}e^{-\sqrt{2\alpha}(\rho-1)a(t)}+e^{-\sqrt{2\alpha}(\rho-1)a(t)}\\
        \le & \left(\sqrt{\frac{\rho+1}{\sqrt{2\alpha}\pi}}\frac{1}{\rho-1}+1\right)e^{-2\alpha(\rho-1)(1- \theta)t}.
	\end{align*}
	Here in the second inequality we used \eqref{6.1},
 \eqref{easy-domi-below}
 and the fact that
 $$
   q\frac{2}{\sqrt{2\alpha}(\rho+1)}=\frac{2\alpha(\rho^2-1)}{\sqrt{2\alpha}(\rho+1)}=\sqrt{2\alpha}(\rho-1).
$$

	(iv) Finally, if  $0< a(t)\le 1$  (that is, $1-\frac{1}{\sqrt{2\alpha}t}\le \theta<1$), then
	$${\bP}\left({M}^Z_t\le \sqrt{2\alpha}\theta t-\sqrt{t}\right)\le 1\le e^{\sqrt{2\alpha}(\rho-1)}e^{-\sqrt{2\alpha}(\rho-1)a(t)}=e^{\sqrt{2\alpha}(\rho-1)}e^{-2\alpha(\rho-1)(1- \theta)t}.$$
	
	The proof is now complete.
	\hfill$\Box$
	
	\medskip

	Recall that $m_t=\sqrt{2\alpha}t-\frac{3}{2\sqrt{2\alpha}}\log t$. The next lemma gives another estimate of $v(t,z)$. The proof  will be given in Appendix.
	
	\begin{lemma}\label{lem_upper-w}
		For any
		$\epsilon\in(0,\sqrt{2\alpha}(\rho-1))$,
		there exist $c_\epsilon>1$ and $T_\epsilon\ge 1$  such that
		$$
		v(t,m_t-z)\le \bP\left(M_t^Z\le m_t-z\right)
		\le c_\epsilon e^{-\sqrt{2\alpha}(\rho-1)z}
		e^{\epsilon z}, \quad t\ge T_\epsilon, z>0.
		$$
	\end{lemma}

	\section{Proofs of the main results}

	Put $\zeta_f(t,x):=\psi'(\lambda^*+u^*_f(t,x))$.  It is clear that $\zeta_f(t,x)\ge \psi'(\lambda^*)=q$.
	\begin{lemma}\label{lemma1}
		For any $f\in\mathcal{B}^+(\R)$,
		\begin{align*}
			U_{1,f}(t,\sqrt{2\alpha}\delta t) \le \left\{
			\begin{array}{ll}
				e^{-qt}, & \hbox{$\delta\ge 0$;} \\
				\frac{1}{2\sqrt{\pi\alpha}|\delta|}
				t^{-1/2}e^{-(q+\alpha\delta^2)t}, & \hbox{$\delta<0$.}
			\end{array}
			\right.
		\end{align*}
	\end{lemma}
	{\bf Proof:} Since $\zeta_f (t,x)\ge \psi'(\lambda^*)=q$,
	by \eqref{eq-v2},
	we have that
	$$U_{1,f}(t,\sqrt{2\alpha}\delta t)=\mE\left(e^{-\int_0^t \zeta_f(t-r,\sqrt{2\alpha}\delta t-B_r)\,\d s};B_t \le \sqrt{2\alpha}\delta t\right)\le e^{-qt}\mP \Big[B_t\le \sqrt{2\alpha}\delta t\Big].$$
	Thus, the desired result follows easily from \eqref{est:bm} with $x=\sqrt{2\alpha}\delta t$.
	\hfill$\Box$
	
	\medskip
	
Note that by the change of variables $s\to t-s$, we have
	$$U_{2,f}(t, x)=\mE \int_0^t e^{-\int_s^t \zeta_f(r,x-B_{t-r})\,dr}\hat{G}_f(s,x-B_{t-s})\,\d s.$$

	\subsection{Proof of Theorem \ref{thm-case1}:  $\delta\in(1-\rho,1)$}
	
It follows from Lemma \ref{fact4} that, to prove  Theorem \ref{thm-case1},
we only need to consider the limiting property of $v_f(t,\sqrt{2\alpha}\delta t)$.
	Note that
	\begin{align}\label{2.7}
		q+\alpha\delta^2-	2\alpha(\rho-1)(1-\delta)=\alpha(\rho-1+\delta)^2,
	\end{align}
	and
	\begin{align}\label{2.8}
		2\alpha(\rho-1)(1-\delta)\le 2\alpha(\rho-1)<\alpha(\rho^2-1)=q,\quad \delta\in[0,1).
	\end{align}
	It follows from Lemma \ref{lemma1} that for any  $\delta\in(1-\rho,1)$,
	$$\lim_{t\to\infty}\frac{e^{2\alpha(\rho-1)(1-\delta)t}}{t^{3(\rho-1)/2}}
U_{1,f}(t,\sqrt{2\alpha}\delta t)=0.
$$
	Thus,
by the decomposition \eqref{decom-v},
to prove the desired result, it suffices to show that
	$$\lim_{t\to\infty}\frac{e^{2\alpha(\rho-1)(1-\delta)t}}{t^{3(\rho-1)/2}}U_{2,f}(t, \sqrt{2\alpha}\delta t)
	=\frac{a_\delta^{3(\rho-1)/2}}{\sqrt{2\alpha}\rho}
	\int_{-\infty}^\infty e^{-\sqrt{2\alpha}(\rho-1)z}A(w_f(z))\d z,$$
	where $a_\delta=1-\frac{1-\delta}{\rho}$ and $A(\lambda)=\frac{1}{\lambda^* }\psi(\lambda)+\psi'(\lambda^*)(1-\lambda/\lambda^*).$

The result above follows from Lemmas \ref{lem-case1} and  \ref{lem-case1-2}
below.
In Lemma \ref{lem-case1-2},
we will show that for  $\delta\in(1-\rho, 1)$, $$\frac{e^{2\alpha(\rho-1)(1-\delta)t}}{t^{3(\rho-1)/2}}\bP\Big(M_t^I\le \sqrt{2\alpha}\delta t,\tau\notin \Big[\frac{1-\delta}{\rho} t-(\log t)\sqrt{t}, \frac{1-\delta}{\rho}  t+(\log t)\sqrt{t}\Big] \Big)\to 0.$$
Thus, on the event
$\left\{M_t^I\le \sqrt{2\alpha}\delta t\right\}$,
with large probability,
the first branching time of the skeleton happens in the interval
$\left[\frac{1-\delta}{\rho} t-(\log t)\sqrt{t}, \frac{1-\delta}{\rho}  t+(\log t)\sqrt{t}\right].$
	
		\medskip
		
	\begin{lemma}\label{lem-case1}
Let $\delta\in(1-\rho, 1)$ and $\mathcal{I}_{t}=[a_\delta t-(\log t)\sqrt{t}, a_\delta t+(\log t)\sqrt{t}]\cap [0,t]$.
Then for any $f\in\mathcal{H}$,
		\begin{align*}
		&\lim_{t\to\infty}\frac{e^{2\alpha(\rho-1)(1-\delta)t}}{t^{3(\rho-1)/2}}
		{\rm E}
		 \int_{\mathcal{I}_t} e^{-\int_s^t \zeta_f(r,\sqrt{2\alpha}\delta t-B_{t-r})\,dr}\hat{G}_f(s,\sqrt{2\alpha}\delta t-B_{t-s})\,\d s\\
		=&\frac{a_\delta^{3(\rho-1)/2}}{\sqrt{2\alpha}\rho}
		\int^\infty_{-\infty}e^{-\sqrt{2\alpha}(\rho-1)z}A(w_f(z))\d z.
	\end{align*}
	\end{lemma}
	{\bf Proof:}
	In this proof, we always assume that $t\ge1$ is large enough such that  $a_\delta t/2\le a_\delta t-(\log t)\sqrt{t}\le a_\delta t+(\log t)\sqrt{t}\le (1+a_\delta )t/2$.
	Since $\psi'$ is increasing and $\psi''$ is decreasing, it follows that, for any $\lambda\ge0$
	$$q=\psi'(\lambda^*)\le \psi'(\lambda^*+\lambda)\le q+\psi''(\lambda^*)\lambda.$$
	Thus we have, for any $s\in \mathcal{I}_t$,
	\begin{align*}
q(t-s)\le \int_s^{t}\zeta_f (r,\sqrt{2\alpha}\delta t-B_{t-r})\,dr&\le q(t-s)+\psi''(\lambda^*)\int_s^t u^*_f(r,\sqrt{2\alpha}\delta t-B_{t-r})\,dr\nonumber\\
&\le q(t-s)+\psi''(\lambda^*)t k(a_\delta t-(\log t)\sqrt{t}).
	\end{align*}
Here  the last inequality follows from Lemma \ref{lemma:u*}(1) and the fact that the function $k$ is decreasing. By  Lemma \ref{lemma:u*}(1), $\sup_{t>1}e^{qt}k(t)<\infty$,
which implies that
$tk(a_\delta t-(\log t)\sqrt{t})\to0$ as $t\to\infty$.
Thus as $t\to\infty$,
\begin{align}\label{3.3.1}
\mE \int_{\mathcal{I}_t} e^{-\int_s^t \zeta_f(r,\sqrt{2\alpha}\delta t-B_{t-r})\,dr}\hat{G}_f(s,\sqrt{2\alpha}\delta t-B_{t-s})\,\d s
\sim \int_{\mathcal{I}_t}e^{-q(t-s)} \mE[\hat{G}_f(s,\sqrt{2\alpha}\delta t-B_{t-s})]\,\d s.
\end{align}
By the change of variables $s=s(u):=a_\delta t+u\sqrt{t}$,
we get that
	\begin{align}\label{1.5}
	  &\int_{\mathcal{I}_t}e^{-q(t-s)} \mE[
	  \hat{G}_f(s,\sqrt{2\alpha}\delta t-B_{t-s})]\,\d s\nonumber\\
=& \int_{\mathcal{I}_t}
e^{-q(t-s)} \mE[\hat{G}_f(s,m_s+(\sqrt{2\alpha}\delta t-m_s-B_{t-s}))]\,\d s\nonumber\\
	  =&
	   \int_{\mathcal{I}_t}e^{-q(t-s)} \d s\int_{\R}\frac{1}{\sqrt{2\pi}(t-s)}e^{-\frac{(z+m_s-\sqrt{2\alpha}\delta t)^2}{2(t-s)}}\hat{G}_f(s,m_s +z)\,\d z\nonumber\\
	   =&
		\sqrt{t}\int^{\log t}_{-\log t}\frac{e^{-q(1-a_\delta )t} e^{q\sqrt{t}u}}{\sqrt{2\pi(t-s(u))}}du
		\int^\infty_{-\infty}e^{-\frac{(m_{s(u)}+z-\sqrt{2\alpha}\delta t)^2}{2(t-s(u))}}\hat{G}_f(s(u),m_{s(u)}+z)\,\d z.
	\end{align}
For $u\in(-\log t ,\log t )$, we have that
	\begin{align*}
		&(m_{s(u)}+z-\sqrt{2\alpha}\delta t)^2  =\left(\sqrt{2\alpha}(a_\delta-\delta) t+\sqrt{2\alpha}u\sqrt{t}-\frac{3}{2\sqrt{2\alpha}}\log (a_\delta t+u\sqrt{t})+z\right)^2 \\
		&=2\alpha(a_\delta-\delta)^2 t^2+2\alpha u^2t+4\alpha(a_\delta-\delta)u t\sqrt{t}-3(a_\delta-\delta)t \log (a_\delta t)\\
		&\quad+2\sqrt{2\alpha}(a_\delta-\delta) zt+R_1(t,u,z),
	\end{align*}
	where $R_1(t,u,z)=\left(-\frac{3}{2\sqrt{2\alpha}}\log (a_\delta t+u\sqrt{t})+z\right)^2-3u\sqrt{t}\log (a_\delta t+u\sqrt{t})+2\sqrt{2\alpha}u\sqrt{t}z- 3(a_\delta-\delta)t \log (1+u/(a_\delta \sqrt{t}))$. Using this one can check that for $|u|\le \log t $,
	\begin{align*}
		R_1(t,u,z)\ge -3(\log t) ^2\sqrt{t}-2\sqrt{2\alpha}(\log t)\sqrt{t}|z|
		-\frac{3(a_\delta-\delta)}{a_\delta} \sqrt{t}\log t.
	\end{align*}
Using the Taylor expansion of $(1-x)^{-1}$, we obtain that
	\begin{align*}
		&\frac{1}{2(t-s(u))}
		 =\frac{1}{2(1-a_\delta)t}\frac{1}{1-u/[(1-a_\delta)\sqrt{t}]}\\
		&=\frac{1}{2(1-a_\delta)t}\left(1+\frac{u}{(1-a_\delta)\sqrt{t}}+\frac{u^2}{(1-a_\delta)^2 t}+R_2(t,u)\right),
	\end{align*}
	where
	$$
	|R_2(t,u)|=\left|\sum_{n=3}^\infty\left[\frac{u}{(1-a_\delta)\sqrt{t}}\right]^n \right|\le \sum_{n=3}^\infty\left[\frac{\log t}{(1-a_\delta)\sqrt{t}}\right]^n
	\le   \frac{2}{(1-a_\delta)^3}(\log t )^3t^{-3/2},
	$$
	here we used the fact that $\log t/[(1-a_\delta)\sqrt{t}]\le 1/2,$ and for $0\le x\le 1/2$, $\sum_{n=3}^\infty x^n=\frac{x^3}{1-x}\le 2x^3$.
	Using the above estimates, we get that for $u\in(-\log t ,\log t )$,
	\begin{align}\label{3.3.2}
		&\frac{(m_{s(u)}+z-\sqrt{2\alpha}\delta t)^2}{2(t-s(u))}\nonumber\\
		=&\frac{\alpha(a_\delta-\delta)^2}{1-a_\delta}\left(t+\frac{u}{(1-a_\delta)}\sqrt{t}+\frac{u^2}{(1-a_\delta)^2 }\right)-\frac{3(a_\delta-\delta)}{2(1-a_\delta)} \log t\nonumber\\
		&+\frac{4\alpha(a_\delta-\delta)u }{2(1-a_\delta)}\left(\sqrt{t}+\frac{u}{(1-a_\delta)}\right)+\frac{2\alpha u^2-3(a_\delta-\delta) \log (a_\delta)+2\sqrt{2\alpha}(a_\delta-\delta) z}{2(1-a_\delta)}\nonumber\\
&+R_3(t,u,z)\nonumber\\
		=&\frac{\alpha(\rho-1)^2(1-\delta)}{\rho}t+qu\sqrt{t}-\frac{3(\rho-1)}{2} \log t+\frac{\alpha\rho^3}{1-\delta}u^2+\sqrt{2\alpha}(\rho-1)z\nonumber\\
		&-\frac{3}{2}(\rho-1)\log(a_\delta)+R_3(t,u,z),
	\end{align}
	where $\lim_{t\to\infty}R_3(t,u,z)=0$ and
	there exists
	a positive function $r(\cdot)$ with $\lim_{t\to\infty}r(t)=0$ such that for any $u\in(-\log t ,\log t )$,
	\begin{align}\label{1.6}
		-R_3(t,u,z)\le r(t)(1+|z|).
	\end{align}
	For any $\epsilon>0$, choose $t_\epsilon$ such that
	$r(t)\le \epsilon$ for any $t>t_\epsilon$.
Noticing that $q(1-a_\delta)+\frac{\alpha(\rho-1)^2(1-\delta)}{\rho}=2\alpha(\rho-1)(1-\delta),$ by \eqref{3.3.1}, \eqref{1.5} and \eqref{3.3.2}, we get that
	\begin{align}\label{1.7}
		&\lim_{t\to\infty}\frac{e^{2\alpha(\rho-1)(1-\delta)t}}{t^{3(\rho-1)/2}}
		\mE \int_{I_t} e^{-\int_s^t \zeta(r,\sqrt{2\alpha}\delta t-B_{t-r})\,dr}\hat{G}_f(s,\sqrt{2\alpha}\delta t-B_{t-s})\,\d s\nonumber\\
		=&a_\delta^{3(\rho-1)/2}\lim_{t\to\infty}\int^{\log t }_{-\log t }\frac{\sqrt{t}}{\sqrt{2\pi(t-s(u))}}e^{-\frac{\alpha\rho^3}{1-\delta}u^2}du
		\int^\infty_{-\infty}e^{-\sqrt{2\alpha}(\rho-1)z}e^{-R_3(t,u,z)}\hat{G}_f(s(u),m_{s(u)}+z)\d z.
	\end{align}
	It follows from \eqref{M-w} that
$$\lim_{t\to\infty}\E\left( e^{-\int_\R f(y-m_t-z)X_t(dy)}; M_t\le m_t+z|\mathcal{S}\right)=\frac{e^{-w_f(z)}-e^{-\lambda^*}}{1-e^{-\lambda^*}}.$$
	Thus by \eqref{e:1}, we get that
	$$\lim_{t\to\infty} v_f(t,m_t+z)=1-\frac{w_f(z)}{\lambda^*}:=\tilde{w}_f(z), $$
	here we used the fact that
	$0\le 1-e^{-u^*(t,x)}\le \P^*(X_t\neq 0)\to 0.$
	It follows that
	\begin{align}\label{2.4}
		&\lim_{t\to\infty}\hat{G}_f(a_\delta t+u\sqrt{t},m_{a_\delta t+u\sqrt{t}}+z)=\lim_{t\to\infty}\hat{G}_f(t,m_t+z)\nonumber\\
		&=
		\frac{1}{\lambda^*}\psi(\lambda^*(1-\tilde{w}_f(z)))+q\tilde{w}_f(z)=A(w_f(z)).
	\end{align}
	Thus, as $t\to\infty$,  the limit of the integrand  in \eqref{1.7} is
	$$
	\frac{a_\delta^{3(\rho-1)/2}}{\sqrt{2\pi(1-a_\delta)}}e^{-\frac{\alpha\rho^3}{1-\delta}u^2-\sqrt{2\alpha}(\rho-1)z}A(w_f(z)).
	$$
	By \eqref{1.6}, \eqref{est-hatG} and Lemma \ref{lem_upper-w},
	we have that, for $\eta$ small enough,
	there exist $T_\eta>1$ and $c_\eta>0$ such that for $t>T_\eta+t_\epsilon$, the integrand in \eqref{1.7}  is smaller than
	$$
	q\frac{a_\delta^{3(\rho-1)/2}}{\sqrt{\pi(1-a_\delta)}}
	e^{-\frac{\alpha\rho^3}{1-\delta}u^2-\sqrt{2\alpha}(\rho-1)z}e^{\epsilon(1+|z|)}\times \left\{
	\begin{array}{ll}
		c_\eta^2 e^{2[\sqrt{2\alpha}(\rho-1)-\eta]z}, & \hbox{$z<0$;} \\
		1, & \hbox{$z>0$,}
	\end{array}
	\right.
	$$
	which is integrable over $\R\times \R$ if we  choose $\epsilon<\sqrt{2\alpha}(\rho-1)$ and $-2\eta+\sqrt{2\alpha}(\rho-1)-\epsilon>0$.
	Thus  using the dominated convergence theorem in \eqref{1.7}, we have that
	\begin{align*}
&\lim_{t\to\infty}\frac{e^{2\alpha(\rho-1)(1-\delta)t}}{t^{3(\rho-1)/2}}
\mE \int_{{\cal I}_t} e^{-\int_s^t \zeta_f(r,\sqrt{2\alpha}\delta t-B_{t-r})\,dr}\hat{G}_f(s,\sqrt{2\alpha}\delta t-B_{t-s})\,\d s\\
=&\frac{a_\delta^{3(\rho-1)/2}}{\sqrt{2\pi(1-a_\delta)}}
\int^{\infty}_{-\infty}e^{-\frac{\alpha\rho^3}{1-\delta}u^2}du\int^\infty_{-\infty}e^{-\sqrt{2\alpha}(\rho-1)z}A(w_f(z))\d z\\
=&\frac{a_\delta^{3(\rho-1)/2}}{\sqrt{2\alpha}\rho}
\int^\infty_{-\infty}e^{-\sqrt{2\alpha}(\rho-1)z}A(w_f(z))\d z.
	\end{align*}
	\hfill$\Box$

	\begin{lemma}\label{lem-case1-2}
		For $\delta\in(1-\rho, 1)$, it holds that
		for any $f\in\mathcal{B}^+(\R)$,
		$$\lim_{t\to\infty}\frac{e^{2\alpha(\rho-1)(1-\delta)t}}{t^{3(\rho-1)/2}}{\rm E}
		\int_{[0,t]\setminus \mathcal{I}_{t}}
		e^{-\int_s^{t}\zeta_f(r,\sqrt{2\alpha}\delta t-B_{t-r}\,dr)}\hat{G}_f(s,\sqrt{2\alpha}\delta t-B_{t-s})\,\d s
		=0.$$
	\end{lemma}
	{\bf Proof:} Since $\zeta_f(t,x)\ge q$,
	using \eqref{est-hatG}
	and the fact that $v_f(t,x)\le v(t,x)$,
	we only need to show that
	\begin{equation}\label{limit-0}
		\lim_{t\to\infty}\frac{e^{2\alpha(\rho-1)(1-\delta)t}}{t^{3(\rho-1)/2}}\mE
		\int_{[0,t]\setminus {\cal I}_{t}}
		e^{-q(t-s)}v^2(s,\sqrt{2\alpha}\delta t-B_{t-s})\,\d s
		=0.
	\end{equation}
Note that
	\begin{align*}
		[0,t]\setminus {\cal I}_t\subset &[0,\epsilon t]\cup
		\left([(a_\delta-\epsilon)t, a_\delta t-(\log t)\sqrt{t}])\cup [a_\delta t+(\log t)\sqrt{t}, (a_\delta+\epsilon)t]\right)\\
		&\cup\left([\epsilon t,(a_\delta-\epsilon)t]\cup[(a_\delta+\epsilon)t,t]\right).
	\end{align*}
	The proof of  \eqref{limit-0} is accomplished  in the following three lemmas by handling the integral over $[0,\epsilon t]$,
	$[(a_\delta-\epsilon)t, a_\delta t-(\log t)\sqrt{t}]\cup [a_\delta t+(\log t)\sqrt{t}, (a_\delta+\epsilon)t]$ and $[\epsilon t,(a_\delta-\epsilon)t]\cup[(a_\delta+\epsilon)t,t]$
	separately. \hfill$\Box$
	
	\begin{lemma}
		Let $\delta\in(1-\rho, 1)$.
		For $\epsilon>0$ small enough,
		$$
		\lim_{t\to\infty}\frac{e^{2\alpha(\rho-1)(1-\delta)t}}{t^{3(\rho-1)/2}}{\rm E}\int_{0}^{\epsilon t} e^{-q(t-s)}v^2(s,\sqrt{2\alpha}\delta t-B_{t-s})\,\d s=0.
		$$
	\end{lemma}
	{\bf Proof:}
	By \eqref{domi-B'},  we have that
	 $$v(s,\sqrt{2\alpha}\delta t-B_{t-s})\le
		 \mP_{B_{t-s}}(B_s\le \sqrt{2\alpha}\delta t)=\mP[B_t\le \sqrt{2\alpha}\delta t|\sigma(B_r: r\le t-s)].
	 $$
Thus it follows that
	\begin{align}\label{3.5.1}
	\mE(v^2\left(s,\sqrt{2\alpha}\delta t-B_{t-s})\right)&\le \mE\left(v(s,\sqrt{2\alpha}\delta t-B_{t-s})\right)\le
		\mP(B_t\le \sqrt{2\alpha}\delta t).
	\end{align}
	Hence, for any $\epsilon>0$,
	\begin{align}\label{2.11}
		&\mE\int_{0}^{\epsilon t} e^{-q(t-s)}v^2(s,\sqrt{2\alpha}\delta t-B_{t-s})\,\d s \le q^{-1}e^{q\epsilon t}e^{-qt}\mP(B_t\le \sqrt{2\alpha}\delta t)\nonumber\\
		&\le  q^{-1}e^{q\epsilon t}\times \left\{\begin{array}{ll}
			e^{-qt}, & \hbox{$\delta\ge 0$;} \\
			\frac{1}{2\sqrt{\pi\alpha}|\delta|}
			t^{-1/2}e^{-(q+\alpha\delta^2)t}, & \hbox{$\delta<0$,}
		\end{array}
		\right.
	\end{align}
	where in the last inequality we used \eqref{est:bm}.
	Using \eqref{2.7} and \eqref{2.8},  we can choose $\epsilon$ small enough so that
	$$
	2\alpha(\rho-1)(1-\delta)+q\epsilon < \left\{\begin{array}{ll}
		q, & \hbox{$\delta\ge 0$;} \\
		q+\alpha\delta^2, & \hbox{$\delta\in(1-\rho,0)$,}
	\end{array}
	\right.
	$$
	which implies the desired result.
	\hfill$\Box$
	
	\begin{lemma}\label{lem:case1-3}
		Let $\delta\in(1-\rho, 1)$.
		For $\epsilon>0$ small enough,
		\begin{align*}
			\lim_{t\to\infty}\frac{e^{2\alpha(\rho-1)(1-\delta)t}}{t^{3(\rho-1)/2}}{\rm E}\left(\int_{(a_\delta-\epsilon)t}^{a_\delta t-
				(\log t)\sqrt{t}}+\int_{a_\delta t+(\log t)\sqrt{t}}^{(a_\delta+\epsilon)t}\right) e^{-q(t-s)}v^2(s,\sqrt{2\alpha}\delta t-B_{t-s})\,\d s=0.
		\end{align*}
	\end{lemma}
	{\bf Proof:}
	Put $S_t:=(a_\delta-\epsilon,a_\delta-(\log t)/\sqrt{t})\cup(a_\delta+(\log t)/\sqrt{t},a_\delta+\epsilon)$. Recall the definition of $m_t$ given by \eqref{def-m_t}.
    By  the change of variables $s=rt$, applying Lemma \ref{lem_upper-w} for $z>0$ and the fact $v\le 1$ for $z\le 0$,
	we get that, for  $\eta$ small enough, there exists $c_\eta\ge 1$ such that for  $t$ large enough,
	\begin{align*}
		&\mE\left(\int_{(a_\delta-\epsilon)t}^{a_\delta t-(\log t)\sqrt{t}}+\int_{a_\delta t+\log t\sqrt{t}}^{(a_\delta+\epsilon)t}\right) e^{-q(t-s)}v^2(s,\sqrt{2\alpha}\delta t-B_{t-s})\,\d s\\
		=&\mE\left(\int_{(a_\delta-\epsilon)t}^{a_\delta t-(\log t)\sqrt{t}}+\int_{a_\delta t+\log t\sqrt{t}}^{(a_\delta+\epsilon)t}\right) e^{-q(t-s)}v^2\left(s, m(s)-(m(s)-\sqrt{2\alpha}\delta t+B_{t-s})\right)\,\d s\\
		\le& c_\eta^2 t\int_{S_t} e^{-q(1-r)t}
		\mE\left[ e^{-2(\sqrt{2\alpha}(\rho-1)-\eta)(m(rt)-\sqrt{2\alpha}\delta t+B_{(1-r)t})}\wedge1\right]\,dr.
	\end{align*}
        We claim that for any $b_1>b_2>0$,
	\begin{equation}\label{domi-e-B1}
		\mE\left(e^{-b_1(b_2+B_1)}\wedge1\right)\le \frac{1}{\sqrt{2\pi}}\left(\frac{1}{b_1-b_2}+\frac{1}{b_2}\right)e^{-b_2^2/2}.
	\end{equation}
Indeed, the left-hand side of \eqref{domi-e-B1} can be written as
\begin{align*}
\mE\left(e^{-b_1(b_2+B_1)};B_1+b_2>0\right)+\mE(B_1+b_2\le 0).
\end{align*}
	By \eqref{est:bm}, we have that
	\begin{equation*}
\mE(B_1+b_2\le 0)=\mE(B_1>b_2)\le \frac{1}{\sqrt{2\pi}}\frac{1}{b_2}e^{-b_2^2/2}.
	\end{equation*}
By the Girsanov theorem, we have
	\begin{align*}
&\mE\left(e^{-b_1(b_2+B_1)};B_1+b_2>0\right)=e^{-b_1b_2} e^{b_1^2/2}\mE (B_1-b_1+b_2>0)\\
\le&  \frac{1}{\sqrt{2\pi}}\frac{1}{b_1-b_2}e^{-b_1b_2} e^{b_1^2/2}e^{-(b_1-b_2)^2/2}
= \frac{1}{\sqrt{2\pi}}\frac{1}{b_1-b_2}e^{-b_2^2/2}.
	\end{align*}
Now \eqref{domi-e-B1} follows immediately.

 We will use  \eqref{domi-e-B1} with
	$b_1=2(\sqrt{2\alpha}(\rho-1)-\eta)\sqrt{(1-r)t}$ and $b_2=\frac{m(rt)-\sqrt{2\alpha}\delta t}{\sqrt{(1-r)t}}.$
	For $\epsilon\in \left(0, \frac{a_\delta-\delta}{2\rho-1}\wedge(1-a_\delta)\right)$,
	we have for any $r\in S_t\subset(a_\delta-\epsilon,a_\delta+\epsilon)$,
	$$
	\frac{\sqrt{2\alpha}(a_\delta+\epsilon-\delta) }{\sqrt{(1-a_\delta-\epsilon)}}\sqrt{t}\ge b_2
	\ge \frac{\sqrt{2\alpha}(a_\delta-\epsilon-\delta) }{\sqrt{(1-a_\delta+\epsilon)}}\sqrt{t} -\frac{\frac{3}{2\sqrt{2\alpha}}}{\sqrt{(1-a_\delta-\epsilon)}}\frac{\log t}{\sqrt{t}},
	$$
	and
	\begin{align*}
		b_1-b_2&\ge 2(\sqrt{2\alpha}(\rho-1)-\eta)\sqrt{(1-a_\delta-\epsilon)t} -\frac{\sqrt{2\alpha}(a_\delta+\epsilon-\delta) }{\sqrt{(1-a_\delta-\epsilon)}}\sqrt{t}\\
		&=\frac{\sqrt{2\alpha}}{\sqrt{1-a_\delta-\epsilon}}\left[2(\rho-1)(1-a_\delta)-(a_\delta-\delta)-(2\rho-1)\epsilon-\frac{2\eta}{\sqrt{2\alpha}}(1-a_\delta-\epsilon)\right]\sqrt{t}\\
		&\ge \frac{\sqrt{2\alpha}}{\sqrt{1-a_\delta-\epsilon}} \left[a_\delta-\delta-(2\rho-1)\epsilon-\frac{2\eta}{\sqrt{2\alpha}}\right]\sqrt{t},
	\end{align*}
where in the final inequality, we used $(\rho-1)(1-a_\delta)=(a_\delta-\delta)$.
	So if we choose
	$\eta\in \left(0, \sqrt{2\alpha}[a_\delta-\delta-(2\rho-1)\epsilon]/2\right)$, and then for $t$ large enough, $b_1>b_2>0$.
	Thus, using \eqref{domi-e-B1}, we have that, for $t$ large enough and
	$r\in S_t$,
	\begin{align}\label{2.12}
		&\mE\left[ e^{-2(\sqrt{2\alpha}(\rho-1)-\epsilon)(m(rt)-\sqrt{2\alpha}\delta t+B_{(1-r)t})}\wedge1\right]
		\le Ct^{-1/2}
		e^{-\frac{(m(rt)-\sqrt{2\alpha}\delta t)^2}{2(1-r)t}}\nonumber\\
		\le & C t^{-1/2}t^{\frac{3(1-\delta)}{2(1-a_\delta-\epsilon)}}
		e^{-\frac{\alpha(r-\delta)^2}{(1-r)}t}.
	\end{align}
	Here in the last inequality we used
the following facts:
$r\le a_\delta+\epsilon<1$ and
	$$e^{-\frac{(m(rt)-\sqrt{2\alpha}\delta t)^2}{2(1-r)t}}\leq
	(rt)^{\frac{3(r-\delta)}{2(1-r)}}e^{-\frac{\alpha(r-\delta)^2}{(1-r)}t}\le t^{\frac{3(1-\delta)}{2(1-a_\delta-\epsilon)}}e^{-\frac{\alpha(r-\delta)^2}{(1-r)}t} .$$
	By Lemma \ref{fact7}, we have that, for $r\in S_t$,
	$$
q(1-r)+\frac{\alpha(r-\delta)^2}{(1-r)}\ge 2\alpha(\rho-1)(1-\delta)+\alpha \rho^2(a_\delta-r)^2\ge 2\alpha(\rho-1)(1-\delta)+\alpha \rho^2\frac{(\log t )^2}{t}.
$$
	Thus,
	there exists $\theta$ such that
\begin{align*}	
&\frac{e^{2\alpha(\rho-1)(1-\delta)t}}{t^{3(\rho-1)/2}}\mE\left(\int_{(a_\delta-\epsilon)t}^{a_\delta t-\log t\sqrt{t}}+\int_{a_\delta t+\log t\sqrt{t}}^{(a_\delta+\epsilon)t}\right) e^{-q(t-s)}v^2(s,\sqrt{2\alpha}\delta t-B_{t-s})\,\d s\\
\le &C t^\theta e^{-\alpha\rho^2(\log t)^2}\to 0,
\quad \mbox{ as } t\to\infty.
\end{align*}
\hfill$\Box$

	\begin{lemma}\label{lem:case1-4}
		Let $\delta\in(1-\rho, 1)$.
		For $\epsilon>0$ small enough,
		\begin{align*}
			\limsup_{t\to\infty}\frac{e^{2\alpha(\rho-1)(1-\delta)t}}{t^{3(\rho-1)/2}}
			{\rm E}\left(\int_{\epsilon t}^{(a_\delta-\epsilon)t}+\int_{(a_\delta+\epsilon)t}^{t}\right)
			e^{-q(t-s)}v^2(s,\sqrt{2\alpha}\delta t-B_{t-s})\,\d s=0.
		\end{align*}
	\end{lemma}
	{\bf Proof:} Set $\mathcal{I}=(\epsilon,a_\delta-\epsilon)\cup (a_\delta+\epsilon,1)$. By the change of variables $r=s/t$, we get that
	\begin{align*}
		&\mE\left(\int_{\epsilon t}^{(a_\delta-\epsilon)t}+\int_{(a_\delta+\epsilon)t}^{t}\right) e^{-q(t-s)}v^2(s,\sqrt{2\alpha}\delta t-B_{t-s})\,\d s\\
		=&t\mE\int_{\mathcal{I}} e^{-q(1-r)t}v^2(rt,\sqrt{2\alpha}\delta t-B_{t-rt})\,\d r\\
		=&t\mE\int_{\mathcal{I}} e^{-q(1-r)t}\, \d r \int_{\R}\frac{1}{\sqrt{2\pi(1-r)t}}e^{-\frac{(z-\sqrt{2\alpha}\delta t)^2}{2(1-r)t}} v^2(rt,z)\,\d z\\
		=&\sqrt{2\alpha}t^2\int_{\mathcal{I}} \,dr\int_{\R}\frac{r}{\sqrt{2\pi(1-r)t}} e^{-q(1-r)t}e^{-\frac{(\sqrt{2\alpha}art-\sqrt{rt}-\sqrt{2\alpha}\delta t )^2}{2(1-r)t}}v^2(rt,\sqrt{2\alpha}\theta rt-\sqrt{rt})\,\d \theta\\
		=&\frac{\sqrt{\alpha}}{\sqrt{\pi}}t^{3/2}
		\int_{\mathcal{I}}\frac{r\,\d r}{\sqrt{1-r}}\left(\int_{-\infty}^{1-\rho}+\int_{1-\rho}^1+\int_{1}^\infty\right)
		e^{-\frac{\alpha\left(\theta r-\frac{\sqrt{r}}{\sqrt{2\alpha t}}-\delta \right)^2t}{(1-r)}-q(1-r)t}v^2(rt,\sqrt{2\alpha}\theta rt-\sqrt{rt})\,\d\theta\\
		=&:\frac{\sqrt{\alpha}}{\sqrt{\pi}}(I_1(t)+I_2(t)+I_3(t)).
	\end{align*}

	For $I_1(t)$, by Lemma \ref{lemma:key}(2)
with $t$ replaced by $rt$,
we have that for $\epsilon t>t_0$ and $\theta<1-\rho$,
	$$v(rt,\sqrt{2\alpha}\theta rt-\sqrt{rt})\le c rt e^{-\alpha \theta^2 rt}e^{-qrt}.$$
	Then by the change of variables $\theta\to -\theta$ in $I_1(t)$, we get that for $t>t_0/\epsilon$,
	\begin{align*}
		I_1(t)&\le
		c^2
		t^{7/2}\int_{\mathcal{I}} \frac{r^3\, dr}{\sqrt{1-r}}\int^{\infty}_{\rho-1} \exp\Big\{-\Big[q(1+r)+\frac{\alpha\left(\theta r+\frac{\sqrt{r}}{\sqrt{2\alpha t}}+\delta \right)^2}{(1-r)}+2\alpha \theta^2 r\Big]t\Big\}\,\d\theta\\
		&\le   c^2
		t^{7/2} e^{-q(1+\epsilon)t}e^{-\alpha\delta^2t}\int_{\mathcal{I}}\frac{r^3\,\d r}{\sqrt{1-r}} \int^{\infty}_{-\infty} e^{-2\alpha \theta^2 rt}\,\d\theta\\
		&= c^2
		t^{7/2} e^{-q(1+\epsilon)t}e^{-\alpha\delta^2t}\int_{\mathcal{I}}\sqrt{\frac{\pi}{2\alpha rt}}\frac{r^3\,dr}{\sqrt{1-r}}\,\d r
		\le Ct^3 e^{-q\epsilon t}
		e^{-(q+\alpha\delta^2)t}.
	\end{align*}
	Since $q+\alpha \delta^2>2\alpha(\rho-1)(1-\delta)$,  it holds that
	$$\lim_{t\to\infty}\frac{e^{2\alpha(\rho-1)(1-\delta)t}}{t^{3(\rho-1)/2}}I_1(t)=0.$$
	
	For $I_2(t)$, by Lemma \ref{lemma:key}(2) and the change of variables $\theta-\frac{1}{\sqrt{2\alpha r t}}\to \theta$, we get
	that for $\epsilon t>t_0$, $I_2(t)$ is less than or equal to
	\begin{align*}
&c^2t^{7/2}
		\int_{\mathcal{I}}  \frac{r^3}{\sqrt{1-r}}\,\d r\int_{1-\rho}^{1} \exp\Big\{-\Big[q(1-r)+\frac{\alpha\left(\theta r-\frac{\sqrt{r}}{\sqrt{2\alpha t}}-\delta \right)^2}{(1-r)}+4\alpha (\rho-1)(1-\theta)r\Big]t\Big\}\,\d\theta\\
	= &c^2t^{7/2}\int_{\mathcal{I}}  \frac{r^3}{\sqrt{1-r}}\,\d r\int_{1-\rho-\frac{1}{\sqrt{2\alpha r t}}}^{1-\frac{1}{\sqrt{2\alpha r t}}} e^{-\left[q(1-r)+\frac{\alpha(\theta r-\delta )^2}{(1-r)}+4\alpha (\rho-1)(1-\theta)r\right]t}e^{2\sqrt{2\alpha}(\rho-1)\sqrt{rt}}\,\d\theta\\
		\le&  Ct^{7/2}e^{2\sqrt{2\alpha}(\rho-1)\sqrt{t}}
		e^{-\inf_{r\in\mathcal{I},\theta<1}H(\theta,r) t},
	\end{align*}
	where $H(\theta,r):=q(1-r)+\frac{\alpha(\theta r-\delta )^2}{(1-r)}+4\alpha (\rho-1)(1-\theta)r$.
	We claim that
	\begin{align}\label{2.10}
		&\inf_{r\in\mathcal{I},\theta<1}H(\theta,r)
		> 2\alpha(\rho-1)(1-\delta).
	\end{align}
	Then it follows  that
	$$\lim_{t\to\infty}\frac{e^{2\alpha(\rho-1)(1-\delta)t}}{t^{3(\rho-1)/2}}I_2(t)=0.$$
	Now we prove \eqref{2.10}. Note that
	\begin{align*}
		H(\theta,r)=&\frac{\alpha r^2}{1-r}\left(\theta-\frac{\delta+2(\rho-1)(1-r)}{r}\right)^2-\alpha (\rho-1)(3\rho-1)(1-r)+4\alpha(\rho-1)(1-\delta).
	\end{align*}
	For $r^*:=\frac{\delta+2(\rho-1)}{2\rho-1}\le r<1$ (that is $\frac{\delta+2(\rho-1)(1-r)}{r}\le 1$) and $\theta<1$, $$H(\theta,r)\ge -\alpha (\rho-1)(3\rho-1)(1-r^*)+4\alpha(\rho-1)(1-\delta) =2\alpha(\rho-1)(1-\delta)+\frac{\alpha(\rho-1)^2(1-\delta)}{2\rho-1} .$$
	For $r\in [0,r^*]\cap \mathcal{I}$ and $\theta<1$, since $\frac{\delta+2(\rho-1)(1-r)}{r}\ge 1$, we have that
	\begin{align*}
		H(\theta,r)&\ge H(1,r)=q(1-r)+\frac{\alpha(r-\delta )^2}{(1-r)}\\
	&\ge 2\alpha(\rho-1)(1-\delta)+\alpha\rho^2(a_\delta-r)^2\\
		&\ge 2\alpha(\rho-1)(1-\delta)+\alpha\rho^2\epsilon^2,
	\end{align*}
where in the second inequality we used Lemma \ref{fact7}.	Thus \eqref{2.10} is valid.

	Finally, we deal with $I_3(t)$.  Since $v(t,x)\le 1$, we have
	\begin{align}\label{3.7.1}
		I_3(t)&\le t^{3/2}
		\int_{\mathcal{I}}\frac{r\,\d r}{\sqrt{1-r}}\int_1^\infty
		e^{-\frac{\alpha\left(\theta r-\frac{\sqrt{r}}{\sqrt{2\alpha t}}-\delta \right)^2t}{(1-r)}-q(1-r)t}\,\d\theta\nonumber\\
		&=  \frac{1}{\sqrt{2\alpha}}t\int_{\mathcal{I}}\,dr\int_{\frac{\sqrt{2\alpha t}(r-\delta)-\sqrt{r}}{\sqrt{1-r}}}^\infty e^{-q(1-r)t} e^{-z^2/2}\,\d z \nonumber\\		
		&\le
		 \frac{\sqrt{\pi}}{\sqrt{\alpha}}t\int_{\mathcal{I}} e^{-q(1-r)t}
\mP\left(B_1\ge \frac{\sqrt{2\alpha t}(r-\delta)-1}{\sqrt{1-r}}\right)\,\d r.
	\end{align}
	If $r\le \delta+\frac{2}{\sqrt{2\alpha t}}$, then
	\begin{align}\label{2.9}
		 &e^{-q(1-r)t} \mP\left(B_1\ge \frac{\sqrt{2\alpha t}(r-\delta)-1}{\sqrt{1-r}}\right)
		\le e^{-q(1-r)t}\le e^{-q(1-\delta)t}e^{\frac{2q}{\sqrt{2\alpha }}\sqrt{t}}\nonumber\\
		= &e^{-2\alpha(\rho-1)(1-\delta)t}e^{-\alpha(\rho-1)^2(1-\delta)t}e^{\frac{2q}{\sqrt{2\alpha }}\sqrt{t}}.
	\end{align}
	If $\delta+\frac{2}{\sqrt{2\alpha t}}< r<1$, then $\frac{\sqrt{2\alpha t}(r-\delta)-1}{\sqrt{1-r}}>1$, and thus by \eqref{est:bm},
	\begin{align}\label{2.9-b}
		e^{-q(1-r)t} \mP\left(B_1\ge \frac{\sqrt{2\alpha t}(r-\delta)-1}{\sqrt{1-r}}\right)
		&\le \frac{1}{\sqrt{2\pi}}\frac{\sqrt{1-r}}{\sqrt{2\alpha t}(r-\delta)-1}e^{-q(1-r)t}e^{-\frac{(\sqrt{2\alpha t}(r-\delta)-1)^2}{2(1-r)}}\nonumber\\
		&\le e^{-q(1-r)t}e^{-\frac{\alpha\left(r-\delta-\frac{1}{\sqrt{2\alpha t}} \right)^2t}{(1-r)}}.
	\end{align}
	It follows from Lemma  \ref{fact7} that for $r\in\mathcal{I}$,
	\begin{align*}
		&q(1-r)+\frac{\alpha(r-\delta -\frac{1}{\sqrt{2\alpha t}})^2}{(1-r)}\nonumber\\
		\ge&2\alpha(\rho-1)\left(1-\delta-\frac{1}{\sqrt{2\alpha t}}\right)+\alpha\rho^2\left(a_\delta-r+\frac{1}{\rho\sqrt{2\alpha t}}\right)^2\nonumber\\
		\ge&2 \alpha(\rho-1)(1-\delta)+\alpha\rho^2\left(\epsilon-\frac{1}{\sqrt{2\alpha t}\rho}\right)^2-\sqrt{2\alpha}(\rho-1)t^{-1/2}.
	\end{align*}
	Then we continue the estimates in \eqref{2.9-b} to get that,
	if $\delta+\frac{2}{\sqrt{2\alpha t}}< r<1$, then
	\begin{equation}\label{2.9-b'}
		e^{-q(1-r)t} \mP\left(B_1\ge \frac{\sqrt{2\alpha t}(r-\delta)-1}{\sqrt{1-r}}\right)\le
		e^{-2\alpha(\rho-1)(1-\delta)t}e^{-\alpha\rho^2\left(\epsilon-\frac{1}{\sqrt{2\alpha t}\rho}\right)^2t+\sqrt{2\alpha}(\rho-1)\sqrt{t}}.
	\end{equation}
	Combining \eqref{3.7.1}, \eqref{2.9} and \eqref{2.9-b'},  we get
	$$
	\limsup_{t\to\infty}\frac{e^{2\alpha(\rho-1)(1-\delta)t}}{t^{3(\rho-1)/2}}I_3(t)=0.
	$$
	The proof is now complete.
	\hfill$\Box$

	\subsection{Proof of Theorem \ref{them:case2}: $\delta=1-\rho$}
	
It follows from Lemma \ref{fact4} that, to prove Theorem \ref{them:case2},
	we only need to consider the limiting property of $v_f(t,\sqrt{2\alpha}\delta t)$.
	It follows from Lemma \ref{lemma1} that for $\delta=1-\rho<0$,
	$$\lim_{t\to\infty}t^{-3(\rho-1)/4}e^{(q+\alpha(1-\rho)^2)t}U_{1,f}(t,\sqrt{2\alpha}(1-\rho) t)=0.$$
Thus,
by the decomposition \eqref{decom-v},
 to prove the desired result, it suffices to show that
\begin{align*}
&\lim_{t\to\infty}t^{-3(\rho-1)/4}e^{(q+\alpha(\rho-1)^2)t}U_{2,f}(t, \sqrt{2\alpha}(\rho-1) t)
\\	=&\frac{1}{\sqrt{2\pi }}\int^{\infty}_{ 0}s^{3(\rho-1)/2}e^{-\alpha\rho^2 s^2}\d s\int^\infty_{-\infty} e^{-\sqrt{2\alpha}(\rho-1)z}A(w_f(z))\d z.
\end{align*}
The display above follows from
 Lemmas \ref{lem-case2} and \ref{lem-case2-2} below.
	In Lemma \ref{lem-case2-2}, we will show that
$$
t^{-3(\rho-1)/4}e^{(q+\alpha(1-\rho)^2) t}\bP\Big(M_t^I\le \sqrt{2\alpha}(1-\rho) t,
\tau\notin \Big[t-(\log t)\sqrt{t},t-t^{1/4}\Big] \Big)\to 0.
$$
	Thus, on the event
$\left\{M_t^I\le \sqrt{2\alpha}(1-\rho) t\right\}$,
with large probability,
the first branching time of the skeleton should happens in  the interval  $\Big[t-(\log t)\sqrt{t},t-t^{1/4}\Big].$

		\medskip
		
	\begin{lemma}\label{lem-case2}
It holds that for any $f\in\mathcal{H}$,
		\begin{align*}
			&\lim_{t\to\infty}t^{-3(\rho-1)/4}e^{(q+\alpha(1-\rho)^2) t}
			{\rm E}
			\int_{t^{1/4}}^{(\log t)\sqrt{t}} e^{-\int_s^t \zeta_f(r,\sqrt{2\alpha}(1-\rho) t-B_{t-r})\,dr}\hat{G}_f(s,\sqrt{2\alpha}(1-\rho) t-B_{t-s})\,\d s\\
			=&\frac{1}{\sqrt{2\pi }}\int^{\infty}_{ 0}s^{3(\rho-1)/2}e^{-\alpha\rho^2 s^2}\d s\int^\infty_{-\infty} e^{-\sqrt{2\alpha}(\rho-1)z}A(w_f(z))\d z.
		\end{align*}
	\end{lemma}
	{\bf Proof:}
In this proof, we always assume that $t\ge1$ is large enough such that  $\log t\le \sqrt{t}$.	Using an argument similar to  that in the first paragraph of the proof of Lemma \ref{lem-case1}, we get that, as $t\to\infty$,
	\begin{align}\label{3.8.1}
		&\mE \int_{t^{1/4}}^{(\log t)\sqrt{t}} e^{-\int_s^t \zeta_f(r,\sqrt{2\alpha}(1-\rho) t-B_{t-r})\,dr}\hat{G}_f\left(s,\sqrt{2\alpha}(1-\rho) t-B_{t-s}\right)\,\d s\nonumber\\
		\sim&\mE \int_{t^{1/4}}^{(\log t)\sqrt{t}} e^{-q(t-s)}\hat{G}_f\left(s,\sqrt{2\alpha}(1-\rho) t-B_{t-s}\right)\,\d s\nonumber\\
		=&\sqrt{t}\int^{\log t}_{t^{-1/4}}\frac{e^{-q(t-u\sqrt{t})}}{\sqrt{2\pi(t-u\sqrt{t})}}\d u\int_{\R}e^{-\frac{\left(m_{u\sqrt{t}}+z+\sqrt{2\alpha}(\rho-1) t\right)^2}{2(t-u\sqrt{t})}}\hat{G}_f\left(u\sqrt{t},m_{u\sqrt{t}}+z\right)\d z.
	\end{align}
For $u\in(t^{-1/4},\log t )$, we have that
\begin{align*}
	(m_{u\sqrt{t}}+z+\sqrt{2\alpha}(\rho-1) t)^2  =&\left(\sqrt{2\alpha}(\rho-1) t+\sqrt{2\alpha}u\sqrt{t}-\frac{3}{2\sqrt{2\alpha}}\log (u\sqrt{t})+z\right)^2 \\
	=&2\alpha(\rho-1)^2 t^2+2\alpha u^2t+4\alpha(\rho-1)u t\sqrt{t}-3(\rho-1)t \log (u\sqrt{t})\\
	&+2\sqrt{2\alpha}(\rho-1) zt+R_4(t,u,z),
\end{align*}
where
\begin{align*}
	R_4(t,u,z)&=\left(-\frac{3}{2\sqrt{2\alpha}}\log (u\sqrt{t})+z\right)^2-3u\sqrt{t}\log (u\sqrt{t})+2\sqrt{2\alpha}u\sqrt{t}z\\
	&\ge -3(\log t) ^2\sqrt{t}-2\sqrt{2\alpha}(\log t)\sqrt{t}|z|.
\end{align*}
Using the Taylor expansion of $(1-x)^{-1}$, we obtain that, for $u\in(t^{-1/4},\log t )$,
\begin{align*}
	\frac{1}{2(t-u\sqrt{t})}
	& =\frac{1}{2t}\frac{1}{1-u/\sqrt{t}}=\frac{1}{2t}\left(1+\frac{u}{\sqrt{t}}+\frac{u^2}{ t}+R_5(t,u)\right),
\end{align*}
where
$
|R_5(t,u)|
\le 2(\log t )^3t^{-3/2}.
$
Thus
	\begin{align*}
		\frac{(m_{u\sqrt{t}}+z+\sqrt{2\alpha}(\rho-1) t)^2}{2(t-u\sqrt{t})}
		=&\alpha(\rho-1)^2t+qu\sqrt{t}-\frac{3(\rho-1)}{4} \log t\\
		&+\alpha\rho^2u^2+\sqrt{2\alpha}(\rho-1)z-\frac{3}{2}(\rho-1)\log(u)+R_6(t,u,z).
	\end{align*}
Here $\lim_{t\to\infty}R_6(t,u,z)=0$ and
there is a positive function $r^*(\cdot)$ with $\lim_{t\to\infty}r^*(t)=0$ such that $-R_6(t,u,z)\le r^*(t)(1+|z|)$ for all $u\in(t^{-1/4},\log t )$.
Now, using \eqref{3.8.1}, we get that
\begin{align*}
&\lim_{t\to\infty}t^{-3(\rho-1)/4}e^{(q+\alpha(1-\rho)^2)t}
\mE \int_{t^{1/4}}^{(\log t)\sqrt{t}} e^{-\int_s^t \zeta_f(r,\sqrt{2\alpha}(1-\rho) t-B_{t-r})\,dr}\hat{G}_f\left(s,\sqrt{2\alpha}(1-\rho) t-B_{t-s}\right)\,\d s\\
=&\lim_{t\to\infty}\int^{\log t}_{t^{-1/4}}\frac{\sqrt{t}}{\sqrt{2\pi(t-u\sqrt{t})}} u^{3(\rho-1)/2}e^{-\alpha\rho^2 u^2}\d u\int_{\R}e^{-\sqrt{2\alpha}(\rho-1)z} e^{-R_6(t,u,z)}\hat{G}_f\left(u\sqrt{t},m_{u\sqrt{t}}+z\right)\d z.
\end{align*}
Using an arguments similar to those in the proof of Lemma \ref{lem-case1}, the desired result follows from the the dominated convergence theorem.
	\hfill$\Box$
	
	\begin{lemma}\label{lem-case2-2}
		It holds that
		for any $f\in\mathcal{B}^+(\R)$,
\begin{align*}
&\lim_{t\to\infty}t^{-3(\rho-1)/4}e^{(q+\alpha(1-\rho)^2) t}
	{\rm E}
	\int_{[0,t]\setminus (t^{1/4},(\log t)\sqrt{t})} e^{-\int_s^t \zeta_f(r,\sqrt{2\alpha}(1-\rho) t-B_{t-r})\,dr}\hat{G}_f(s,\sqrt{2\alpha}(1-\rho) t-B_{t-s})\,\d s\\
&=0.
\end{align*}
	\end{lemma}
	{\bf Proof:} We only need to show that
		\begin{align*}
		\lim_{t\to\infty}\frac{e^{(q+\alpha(1-\rho)^2) t}}{t^{3(\rho-1)/4}}
		\mE \int _{(0,t)\setminus (t^{1/4},(\log t)\sqrt{t})}
		e^{-q(t-s)}v^2(s,\sqrt{2\alpha}(1-\rho) t-B_{t-s})\,\d s=0.
	\end{align*}
We prove the above result in three steps.

 {\it Step 1}:
	By \eqref{3.5.1}, we have that
	$$\mE\left(v^2(s,\sqrt{2\alpha}(1-\rho) t-B_{t-s})\right)\le \mP\left(B_t\le \sqrt{2\alpha}(1-\rho) t\right)\le \frac{1}{2\sqrt{\pi\alpha}(\rho-1)\sqrt{t}}e^{-\alpha(\rho-1)^2t}.$$
	Thus,  for any $T>0$,
	\begin{align}\label{3.4}
		&	\frac{e^{(q+\alpha(\rho-1)^2) t}}{t^{3(\rho-1)/4}}\mE\int_0^T e^{-q(t-s)}v^2(s,\sqrt{2\alpha}(1-\rho) t-B_{t-s})\,\d s\nonumber\\
	\le &\int_0^T e^{qs}\,\d s\frac{1}{2\sqrt{\pi\alpha}(\rho-1)\sqrt{t}} \frac{1}{t^{3(\rho-1)/4}}\to 0,
		\quad\mbox{ as } t\to\infty.
	\end{align}

{\it Step 2:}
	Using arguments similar to those in the proofs of Lemmas \ref{lem:case1-3} and \ref{lem:case1-4},
	we get that,
	$$\frac{e^{(q+\alpha(\rho-1)^2)t}}{t^{3(\rho-1)/4}}\mE\int_{\sqrt{t}\log t}^{t} e^{-q(t-s)}v^2(s,\sqrt{2\alpha}(1-\rho) t-B_{t-s})\,\d s\to 0,\quad \mbox{ as }t\to\infty.$$

{\it Step 3:}
	Note that there exists $T_0$ such that $m_s>0$ for all $s>T_0$.
	Using Lemma \ref{lem_upper-w}, we get that, for  $\eta$ small enough, there exist $c_\eta>1$ and $T_\eta>1$  such that for $T>T_\eta+T_0$,
	\begin{align}\label{3.9.2}
		&\mE \int_{T}^{t^{1/4}}e^{-q(t-s)}v^2(s,\sqrt{2\alpha}(1-\rho) t-B_{t-s})\,\d s\nonumber\\
		=&\mE\int_{T}^{t^{1/4}} e^{-q(t-s)}v^2(s, m(s)-(m(s)+\sqrt{2\alpha}(\rho-1) t+B_{t-s}))\,\d s\nonumber\\
		\le &c_\eta^2 \int_{T}^{t^{1/4}} e^{-q(t-s)}
		\mE[ e^{-2(\sqrt{2\alpha}(\rho-1)-\eta)(m(s)+\sqrt{2\alpha}(\rho-1) t+B_{t-s})}\wedge1]\,\d s.
	\end{align}
	Similar to \eqref{2.12}, we have that, for $T<s<t^{1/4}$,
	\begin{align}\label{3.9.1}
		&\mE[ e^{-2(\sqrt{2\alpha}(\rho-1)-\eta)(m(s)+\sqrt{2\alpha}(\rho-1) t+B_{t-s})}\wedge1]\nonumber\\
		\le& Ct^{-1/2}e^{-\frac{(m(s)+\sqrt{2\alpha}(\rho-1) t)^2}{2(t-s)}}\nonumber\\
		\le& Ct^{-1/2}t^{\frac{3(\rho-1)}{8}}e^{-\alpha(\rho-1)^2 t-qs}
	\end{align}
with $C$ being a positive constant.
	Here in the last inequality, we used the fact that
	\begin{align*}
		\frac{(m(s)+\sqrt{2\alpha}(\rho-1) t)^2}{2(t-s)}&= \frac{(\sqrt{2\alpha}\rho s-\frac{3}{2\sqrt{2\alpha}}\log s+\sqrt{2\alpha}(\rho-1) (t-s))^2}{2(t-s)}\nonumber\\
		&\ge \alpha(\rho-1)^2(t-s)+\sqrt{2\alpha}(\rho-1) \left(\sqrt{2\alpha}\rho s-\frac{3}{2\sqrt{2\alpha}}\log s\right)\nonumber\\
		&=\alpha(\rho-1)^2 t+qs-\frac{3}{2}(\rho-1)\log s.
	\end{align*}
Putting \eqref{3.9.1} back to \eqref{3.9.2}, we get that
	\begin{align*}
		&\frac{e^{(q+\alpha(\rho-1)^2) t}}{t^{3(\rho-1)/4}}\mE\int_T^{t^{1/4}} e^{-q(t-s)}v^2(s,\sqrt{2\alpha}(1-\rho) t-B_{t-s})\,\d s
		\le Ct^{-1/4}t^{\frac{-3(\rho-1)}{8}}\to 0,\quad \mbox{as } t\to\infty.
	\end{align*}
	Now the proof is complete.
	\hfill$\Box$
	
	\subsection{ Proof of Theorem \ref{them:case3} : $\delta<1-\rho$}
By \eqref{eq-v2''} we have that
	\begin{align}\label{eq-v}
		v_f(t,x)&=e^{-qt}\mE \Big[B_t\le x\Big]+\mE \int_0^t e^{-q(t-s)}G_f(s,x-B_{t-s})\,\d s,
	\end{align}
	where
	\begin{align}\label{def:G}
		G_f(t,x):&=\hat{G}_f(t,x)-\phi(u^*_f(t,x))v_f(t,x)\\
		&=\frac{1}{\lambda^* }\Big[\psi(\lambda^*+u^*_f(t,x)-\lambda^*v_f(t,x))-\psi(\lambda^*+u^*_f(t,x))\Big]+qv_f(t,x)\nonumber
	\end{align}
	with
	$\phi(\lambda)=\psi'(\lambda+\lambda^*)-q$
	being defined by \eqref{phi}.
	
	\medskip

	It follows from Lemma \ref{fact4} that, to prove Theorem \ref{them:case3},
	we only need to consider the limiting property of $v_f(t,\sqrt{2\alpha}\delta t)$.
Using L'Hospital's rule, one has that
\begin{equation}\label{tail-normal}
\lim_{x\to\infty}\frac{\mP(B_1>x)}{x^{-1}e^{-x^2/2}}=\frac{1}{\sqrt{2\pi}}\lim_{x\to\infty}\frac{\int_x^\infty e^{-y^2/2}\,\d y}{x^{-1}e^{-x^2/2}}=\frac{1}{\sqrt{2\pi}}.
\end{equation}
It follows that
	\begin{align}
		\lim_{t\to\infty}\sqrt{t}e^{(q+\alpha\delta^2)t}	e^{-qt}\mE \Big[B_t\le \sqrt{2\alpha}\delta t\Big]=&\lim_{t\to\infty}\sqrt{t}e^{(q+\alpha\delta^2)t}	e^{-qt}\mE \Big[B_1\ge \sqrt{2\alpha}|\delta| \sqrt{t}\Big]
=\frac{1}{2\sqrt{\pi\alpha}|\delta|}.
	\end{align}
Hence, by  \eqref{eq-v}, to prove the desired result, we only need to prove that
\begin{align*}
&\lim_{t\to\infty}\sqrt{t}e^{(q+\alpha\delta^2)t}\mE\int_0^{t} e^{-q(t-s)}
	G_f(s,\sqrt{2\alpha}\delta t-B_{t-s})\,\d s\\
	=&\frac{1}{\sqrt{2\pi}}\int_0^\infty e^{(q-\alpha\delta^2)s}\,\d s\int_{\R}e^{\sqrt{2\alpha}\delta z}{G}_f(s,z)\,\d z,
\end{align*}
which will follow from Lemmas \ref{int-upto-(t-a)} and  \ref{int-(t-a,t)} below. In Lemma \ref{int-(t-a,t)}, we will show that,
 for any $T>0$,
$$\sqrt{t}e^{(q+\alpha\delta^2)t}\bP\Big(M_t^I\le \sqrt{2\alpha}\delta t,\tau\in [0,t-T] \Big)\to 0.$$
Thus, on the event
$\left\{M_t^I\le \sqrt{2\alpha}\delta t\right\}$,
with large probability,
 the first branching of the skeleton happens in the interval  $[t-T,t]. $

\medskip

	\begin{lemma}\label{int-upto-(t-a)}
If $\delta<1-\rho$, then for any $f\in\mathcal{B}_b^+(\R)$ and any $T>0$, it holds that
\begin{align*}
&\lim_{t\to\infty}\sqrt{t}e^{(q+\alpha\delta^2)t}{\rm  E}\int_0^{t-T} e^{-q(t-s)}{G}_f(s,\sqrt{2\alpha}\delta t-B_{t-s})\,\d s\\
=&\frac{1}{\sqrt{2\pi}}\int_0^\infty e^{(q-\alpha\delta^2)s}\,\d s\int_{\R}e^{\sqrt{2\alpha}\delta z}{G}_f(s,z)\,\d z.
\end{align*}
		
	\end{lemma}
	{\bf Proof:}
	Note that
	\begin{align*}
		&\sqrt{t}e^{(q+\alpha\delta^2)t}
		\mE\int_0^{t-T} e^{-q(t-s)}{G}_f(s,\sqrt{2\alpha}\delta t-B_{t-s})\,\d s\\
		=&\int_0^{t-T} \frac{\sqrt{t}}{\sqrt{2\pi(t-s)}}e^{(q-\alpha\delta^2)s}\,\d s\int_{\R}e^{\sqrt{2\alpha}\delta z}e^{-\frac{(z-\sqrt{2\alpha}\delta s)^2}{2(t-s)}}{G}_f(s,z)\,\d z.
	\end{align*}
	The absolute value of the integrand above is less than
       $\frac{1}{\sqrt{2\pi}}\sqrt{1+s/T}e^{(q-\alpha\delta^2)s}e^{\sqrt{2\alpha}\delta z}|{G}_f(s,z)|$,
	thus by the dominated convergence theorem, it suffices to show that
	\begin{align}\label{4.1}
		\int_{0}^{\infty} \sqrt{s+T}e^{(q-\alpha\delta^2)s}\,\d s\int_\R e^{\sqrt{2\alpha}\delta z}|{G}_f(s,z)|\,\d z<\infty.
	\end{align}
	By \eqref{def:G}, \eqref{est-hatG} and the fact that $v_f(t,x)\le v(t,x)$, we have that
	\begin{equation}\label{est:G}
	|G_f(s,z)|\le \phi(u^*_f(s,z))v_f(s,z)+ \hat{G}_f(s,z)\le \phi(u^*_f(s,z))v(s,z)+qv(s,z)^2.
	\end{equation}
We will prove  \eqref{4.1} in two steps.
Recall that
$k(t)=-\log \P^*(\|X_t\|=0)$.

{\it Step 1}:
First we consider the integral over $s\in(0,A)$, where $A>0$ is a constant.
Since $\phi$ is increasing, by Lemma \ref{lemma:u*}(1), $\phi(u^*_f(s,z))\le \phi(k(s))$. By lemma \ref{lemma:key}(1), $v(s,z)\le\mP(B_s\le z)=\mP(B_1\le z/\sqrt{s})$.
Thus we have for $0<s<A$,
\begin{align*}
&\int_{-\infty}^{0}e^{\sqrt{2\alpha}\delta z}\phi(u^*_f(s,z))v(s,z)\,\d z
 \le \phi(k(s))\int_{-\infty}^{0}e^{\sqrt{2\alpha}\delta z}\mP(B_1\le z/\sqrt{s})\d z\\
 =&\sqrt{s}\phi(k(s))\int_0^{\infty}e^{\sqrt{2\alpha}|\delta| \sqrt{s}z}\mP(B_1\ge  z)\d z
 \le\sqrt{s}\phi(k(s))\int_0^{\infty}e^{\sqrt{2\alpha}|\delta| \sqrt{A}z}\mP(B_1\ge  z)\d z.
\end{align*}
Since $\mP(B_1\ge z)\sim \frac{1}{\sqrt{2\pi}}z^{-1}e^{-z^2/2}$ as $z\to\infty$, we have
$\int_0^{\infty}e^{\sqrt{2\alpha}|\delta| \sqrt{A}z}\mP(B_1\ge  z)\d z <\infty$.
Thus
\begin{align}\label{3.11.2}
	\int_{-\infty}^{0}e^{\sqrt{2\alpha}\delta z}\phi(u^*_f(s,z))v(s,z)\,\d z \le C\sqrt{s}\phi(k(s)).
\end{align}
For any $\epsilon>0$, since $v(s,z)\le 1$, we have
\begin{align}\label{3.11.3}
	\int_{0}^{s^\epsilon}e^{\sqrt{2\alpha}\delta z}\phi(u^*_f(s,z))v(s,z)\,\d z \le s^\epsilon\phi(k(s)).
\end{align}
By \eqref{3.11.2}, \eqref{3.11.3} and Lemma \ref{lemma2}, for any $\epsilon>0$,
\begin{align}\label{3.11.4}
\int_0^A  \sqrt{s+T}e^{(q-\alpha\delta^2)s} \int_{-\infty}^{s^\epsilon}e^{\sqrt{2\alpha}\delta z}\phi(u^*_f (s,z))v(s,z)\,\d z\d s<\infty.
\end{align}
Since $\phi'(\lambda)=\psi''(\lambda^*+\lambda)$ is decreasing and $\phi(0)=0$, we have
\begin{equation}\label{domi-phi}
\phi(\lambda)\le \phi'(0)\lambda.
\end{equation}
Thus, by \eqref{est:u*},
$$\phi(u^*_f(s,z))\le \phi'(0)u^*_f(s,z)\le C (1+z^{-2/\vartheta})e^{(a+\alpha)s},\quad z>0.$$
Since  $v(s,z)\le 1$, we have  for $0<s<A$,
	\begin{align*}
		&\int_{s^\epsilon}^\infty e^{\sqrt{2\alpha}\delta z}\phi(u^*_f(s,z))v(s,z)\,\d z\le  Ce^{(a+\alpha)s} \int_{s^\epsilon}^\infty e^{-\sqrt{2\alpha}|\delta| z} (1+z^{-2/\vartheta})\d z\\
		&\le Ce^{(a+\alpha)A} \Big[\int_{s^\epsilon}^{A^\epsilon }  (1+z^{-2/\vartheta})\d z+ \int_{A^\epsilon}^\infty e^{-\sqrt{2\alpha}|\delta| z} (1+z^{-2/\vartheta})\d z\Big]
		\le
		C(1+s^{\epsilon(1-2/\vartheta)}).
		\end{align*}
	Now we choose $\epsilon$ small enough such that $\epsilon(2/\vartheta-1)<1$. Thus
	\begin{align}\label{3.11.5}
		\int_0^A  \sqrt{s+T}e^{(q-\alpha\delta^2)s} \int_{s^\epsilon}^\infty e^{\sqrt{2\alpha}\delta z}\phi(u^*_f(s,z))v(s,z)\,\d z\d s<\infty.
	\end{align}
Combining \eqref{3.11.4} and \eqref{3.11.5},  we obtain that
	\begin{align}\label{3.11.6}
	\int_0^A  \sqrt{s+T}e^{(q-\alpha\delta^2)s} \int_{-\infty}^\infty e^{\sqrt{2\alpha}\delta z}\phi(u^*_f(s,z))v(s,z)\,\d z\d s<\infty.
\end{align}

{\it Step 2}: By Lemma \ref{lemma:u*}(1), $\sup_{s>A}e^{qs}k(s)=e^{qA}k(A)<\infty$. Hence  we have  for $s>A$,
$$\phi(u^*_f(s,z))\le \phi'(0)u^*_f(s,z)\le \phi'(0)k(s)\le \phi'(0)e^{qA}k(A)e^{-qs}.$$
	Thus we get that, for $s>A$,
	\begin{align}\label{3.11.7}
		&\int_\R e^{\sqrt{2\alpha}\delta z}\phi(u^*_f(s,z))v(s,z)\,\d z\le C e^{-qs}\int_\R e^{\sqrt{2\alpha}\delta z}v(s,z)\,\d z\nonumber\\
		= &C\sqrt{2\alpha} s e^{-qs} e^{-\sqrt{2\alpha}\delta \sqrt{s}}\int_\R e^{2\alpha \delta s \theta }v(s,\sqrt{2\alpha}\theta s-\sqrt{s})\,\d \theta.
	\end{align}
	We will
	divide the above integral into three parts: $\int_{1}^\infty+\int_{1-\rho}^1+\int_{-\infty}^{1-\rho}$. We deal with them one by one. Using Lemma  \ref{lemma:key}(2), we have that for $A>t_0$ and $s>A,$
		 $$\int_1^\infty  e^{2\alpha \delta s \theta }v(s,\sqrt{2\alpha}\theta s-\sqrt{s})\,\d \theta
		 \le \int_1^\infty  e^{-2\alpha |\delta| s\theta}\,\d \theta
		 =\frac{1}{\sqrt{2\alpha}|\delta|s}e^{-{2\alpha}|\delta| s},$$
		\begin{align*}
			\int_{1-\rho}^1  e^{2\alpha \delta s \theta}v(s,\sqrt{2\alpha}\theta s-\sqrt{s})\,\d \theta
			&\le cs \int_{1-\rho}^1  e^{2\alpha \delta s \theta }e^{-2\alpha(\rho-1)(1-\theta)s}\,\d \theta\\
			&\le cs \rho e^{-2\alpha(\rho-1)(\rho+\delta)s},
			\end{align*}
       and
		\begin{align*}
				&\int_{-\infty}^{1-\rho}  e^{2\alpha \delta s \theta }v(s,\sqrt{2\alpha}\theta s-\sqrt{s})\,\d \theta
			\le cs \int_{-\infty}^{1-\rho}   e^{2\alpha \delta s \theta }e^{-(q+\alpha \theta^2)s}\,\d \theta\\
           =&cs e^{(-q+\alpha\delta^2)s}\int_{-\infty}^{1-\rho}   e^{-\alpha s(\theta-\delta)^2}\d \theta
			\le Cs^{1/2}e^{(-q+\alpha\delta^2)s}.
		\end{align*}
	For $\delta<1-\rho$, one can check that $$2\alpha\delta\le -2\alpha (\rho-1)(\rho+\delta)\le -q+{\alpha}\delta^2.$$
	Thus for $s>A$,
	\begin{align}\label{3.11.8}
\int_{-\infty}^{\infty}  e^{2\alpha \delta s \theta }v(s,\sqrt{2\alpha}a s-\sqrt{s})\,\d \theta
\le Cs e^{(-q+\alpha\delta^2)s}.
	\end{align}
	It follows from \eqref{3.11.7}  and \eqref{3.11.8} that
	\begin{align*}
	&	\int_{A}^{\infty} \sqrt{s+T}e^{(q-\alpha\delta^2)s}\,\d s\int_\R e^{\sqrt{2\alpha}\delta z}\phi(u^*(s,z))v(s,z)\,\d z\\
  \le &C \int_{A}^{\infty} \sqrt{s+T}s^{2}e^{-qs}e^{-\sqrt{2\alpha}\delta\sqrt{s}}\,\d s <\infty.
	\end{align*}
	Combining the two steps above, we get
	\begin{align*}
		\int_{0}^{\infty} \sqrt{s+T}e^{(q-\alpha\delta^2)s}\,\d s\int_\R e^{\sqrt{2\alpha}\delta z}\phi(u^*_f(s,z))v(s,z)\,\d z & <\infty.
	\end{align*}
	Similarly, one can prove that
	\begin{align*}
		\int_{0}^{\infty} \sqrt{s+T}e^{(q-\alpha\delta^2)s}\,\d s\int_\R e^{\sqrt{2\alpha}\delta z}v(s,z)^2\,\d z & <\infty.
	\end{align*}
	Hence \eqref{4.1} holds and the desired result follows immediately.
	\hfill$\Box$
	
	\begin{lemma}\label{int-(t-a,t)}
If $\delta<1-\rho$, then for any $f\in\mathcal{B}_b^+(\R)$ and $T>0$,
		$$
		\lim_{t\to\infty}\sqrt{t}e^{(q+\alpha\delta^2)t}{\rm E}\int_{t-T}^t e^{-q(t-s)}G_f(s,\sqrt{2\alpha}\delta t-B_{t-s})\,\d s=0.
		$$
	\end{lemma}
{\bf Proof:}
Note that
\begin{align*}
	&\mE\int_{t-T}^t e^{-q(t-s)}|G_f(s,\sqrt{2\alpha}\delta t-B_{t-s})|\,\d s=\int_{0}^T e^{-qs}\mE|G_f(t-s,\sqrt{2\alpha}\delta t-B_{s})|\,\d s\\
	=&\int_{0}^T e^{-qs}\mE[|G_f(t-s,\sqrt{2\alpha}\delta t-B_{s})|;B_s<-(\epsilon t-\sqrt{t})]\,\d s\\
	&+
	\int_{0}^T e^{-qs}\mE[|G_f(t-s,\sqrt{2\alpha}\delta t-B_{s})|;B_s\ge-(\epsilon t-\sqrt{t})]\,\d s,
\end{align*}
where $\epsilon<1-\rho-\delta$ is a small constant.

By \eqref{domi-phi} and Lemma \ref{lemma:u*}(1), $\sup_{t>1}\phi(u^*_f(t,x))\le \phi'(0)\sup_{t>1}u^*_f(t,x)\le \phi'(0)k(1)<\infty.$  Since $v(t,x)\le 1$, we have
$\sup_{t>1}\sup_{x}|G_f(t,x)|<+\infty$.  Hence we have, for $t>1$ large enough, and $s\in(0,T)$,
\begin{align*}
&\mE\left[|G_f(t-s,\sqrt{2\alpha}\delta t-B_{s})|;B_s\le -(\epsilon t-\sqrt{t})\right]\le C\mP\left(B_s\ge(\epsilon t-\sqrt{t})\right)\\
\le& C\frac{\sqrt{s}}{\epsilon t-\sqrt{t}}e^{-(\epsilon t-\sqrt{t})^2/(2s)}\le C\frac{\sqrt{T}}{\epsilon t-\sqrt{t}}e^{-(\epsilon t-\sqrt{t})^2/(2T)},
\end{align*}
where in the second inequality, we used \eqref{est:bm}.

Thus for any $\epsilon>0$,
as $t\to\infty$,
$$\sqrt{t}e^{(q+\alpha\delta^2)t}\int_{0}^T e^{-qs}\mE\left[|G_f(t-s,\sqrt{2\alpha}\delta t-B_{s})|;B_s<-(\epsilon t-\sqrt{t})\right]\,\d s\to 0.$$
Note that  if $B_s\ge-(\epsilon t-\sqrt{t})$,  then
$$\sqrt{2\alpha}\delta t-B_s\le \sqrt{2\alpha}(\delta+\epsilon) t-\sqrt{t}\le \sqrt{2\alpha}(\delta+\epsilon) (t-s)-\sqrt{t-s}.$$
Using Lemma \ref{lemma:key}(2) with $\theta=\delta+\epsilon<1-\rho$,
for $t>t_0+T$ and $s\in(0,T)$,
$$
v(t-s,\sqrt{2\alpha}\delta t-B_s)\le v(t-s,\sqrt{2\alpha}(\delta+\epsilon)(t-s)-\sqrt{t-s} )\le  c te^{-q(t-s)}e^{-\alpha(\delta+\epsilon)^2(t-s)}.
$$
By Lemma \ref{lemma:u*}(1), we have that for
$t\ge t_0+T$ and $s\in(0,T)$,
\begin{align*}
	&\phi(u^*_f(t-s,\sqrt{2\alpha}\delta t+z))\le \phi'(0)u^*_ft-s,\sqrt{2\alpha}\delta t+z)\\
	\le&\phi'(0) k(t-s)\le \phi'(0)e^{qt_0}k(t_0)e^{-q(t-s)}.
\end{align*}
Thus, by \eqref{est:G}, we get that, if $B_s\ge-(\epsilon t-\sqrt{t})$,
\begin{align}
	&|G_f(t-s,\sqrt{2\alpha}\delta t-B_s)|\le Ct^2e^{-2q(t-s)}e^{-\alpha(\delta+\epsilon)^2(t-s)}\nonumber\\
	\le &Ce^{2qs}e^{\alpha(\delta+\epsilon)^2 s}t^2e^{-2qt}e^{-\alpha\delta^2 t}e^{-2\alpha\delta\epsilon t}.
\end{align}
It follows that, as $t\to\infty$,
\begin{align}
	&\sqrt{t}e^{(q+\alpha\delta^2)t}\int_{0}^T e^{-qs}\mE[|G_f(t-s,\sqrt{2\alpha}\delta t-B_{s})|;B_s\ge-(\epsilon t-\sqrt{t})]\,\d s\nonumber\\
	\le& Ct^{5/2} e^{-(q+2\alpha\delta\epsilon)t}\int_0^T e^{qs}e^{\alpha(\delta+\epsilon)^2 s}\d s\le  Ct^{5/2} e^{-(q+2\alpha\delta\epsilon)t}\to 0,
\end{align}
if we choose $\epsilon$ small enough such that $q+2\alpha\delta\epsilon>0$.
The proof is now complete.
\hfill$\Box$

\appendix
	\section{Appendix}
\begin{lemma}\label{Ztlek}
		For $k\ge1$,
		$$\bP(\|Z_t\|\le k)\le k e^{-qt}.$$
	\end{lemma}
	{\bf Proof:}
	Let $Z'_t$ be a continuous time branching process with branching rate $q$, and when a particle dies, it splits into two particles.
	Then $Z'_t$ is a pure birth process, and the distribution of $Z'_t$ is given by
	$$\bP(Z'_t\le k)=1-(1-e^{-qt})^k.$$
	
According to the definition of  $Z_t$,
each particle  splits into  at least two children ($p_0=p_1=0$),  then we get that
	$$\bP(\|Z_t\|\le k)\le \bP(Z'_t\le k)=1-(1-e^{-qt})^k\le ke^{-qt}.$$
	\hfill$\Box$

	{\bf Proof of Lemma \ref{lem_upper-w}:}
	Since $v(t, m_t-z)\le 1$, it is clear that the desired result is valid  for $z<1$. In the following, we only need to consider the case $z\ge 1.$
	Put $a^*=\sqrt{2\alpha}(\rho-1)/q$. Assume that $\eta\in(0,a^*/2)$ and $t\ge 1$.
	
	(i) First we deal with the case $z>\frac{a^*}{\eta}\sqrt{t}$.
		Since for any $\theta$,
	\begin{equation*}
	q+\alpha \theta^2-	2\alpha(\rho-1)(1-\theta)=\alpha(\rho-1+\theta)^2\ge0,
	\end{equation*}
then by Lemma \ref{lemma:key}(2), one has that there exits $t_0>1$ and $c>0$ such that,  for any $t>t_0$ and $\theta<1$,
\begin{equation}
	v(t,\sqrt{2\alpha}\theta t-\sqrt{t})\le \bP(M_t^Z\le \sqrt{2\alpha}\theta t-\sqrt{t})\le c t e^{-	2\alpha(\rho-1)(1-\theta)t}.
\end{equation}
Thus, using the above inequality with $\theta=1-\frac{z-\sqrt{t}}{\sqrt{2\alpha}t}<1$, we get that for any $t>t_0$,
	\begin{align}\label{1.4}
		\bP(M_t^Z\le m_t-z)&\le  \bP(M_t^Z\le \sqrt{2\alpha}t-z)=\bP(M_t^Z\le \sqrt{2\alpha}\theta t-\sqrt{t})\nonumber\\
		&\le cte^{-\sqrt{2\alpha}(\rho-1)(z-\sqrt{t})}
		\le cz^2e^{-\sqrt{2\alpha}(\rho-1)z}e^{q\eta z},
	\end{align}
	where in the final inequality, we use the fact that
	$t\le (\frac{\eta}{a^*}z)^2\le z^2$
	and $\sqrt{2\alpha}(\rho-1)\sqrt{t}=qa^*\sqrt{t}\le q\eta z.$
	
	(ii) Now we consider the case $z\in[1,\frac{a^*}{\eta}\sqrt{t}]$.    Put
 $K:=[a^*/\eta]$.
Note that $K\ge 1$.
	Define $s_n= \eta  z n$.
	In the following, we always assume that $t$ is large enough such that
	$s_K<t.$
	Note that
	\begin{align}\label{2.2}
		\bP(M_t^Z\le m_t-z)\le
\bP(\|Z_{s_K}\|\le z^2)+ \sum_{l=1}^{K}\bP\left(\|Z_{s_{l-1}}\|\le  z^2<\|Z_{s_{l}}\|,M_t^Z\le m_t-z \right).
	\end{align}
	By Lemma \ref{Ztlek}, we have that
	\begin{equation}\label{1.1}
\bP(\|Z_{s_K}\|\le z^2)\le  z^2 e^{-qs_K}=z^2 e^{-q\eta K z}\le z^2e^{q\eta  z} e^{-\sqrt{2\alpha}(\rho-1)z}.
	\end{equation}
	Now we deal with the second part of the right-hand side of \eqref{2.2}.
Suppose $1\le l\le K$.
	Note that for any  $u\in\mathcal{L}_{s_l},$ $z_u(s_l)\overset{d}{=}Y\sim N(0,s_l).$
 Let $M_{t}^{Z,u}:=\max_{v\in \mathcal{L}_{t},u\preccurlyeq v} z_v(t)-z_u(s_l)$, 	for any $u\in\mathcal{L}_{s_l}$.  By the branching property of $Z$, given  $\sigma(\|Z_s\|,s\in[0,s_l])$, $\{M_{t}^{Z,u}, u\in\mathcal{L}_{s_l}\}$ are i.i.d. with the same distribution as $(M_{t-s_l}^Z,\bP)$, and independent of $\{z_u(s_l), u\in\mathcal{L}_{s_l}\}$.
It is clear that
	$$M_t^Z=\max_{u\in\cL_{s_l}} [z_u(s_l)+M_{t}^{Z,u}].$$
	It follows from \cite[Lemma 5.1]{GH} that
	\begin{align*}
		&	\bP\left(M_t^Z\le m_t-z\,\big|\sigma(\|Z_s\|,s\in[0,s_l])\right)\\
		\le& \bP\left(Y+\max_{u\in\cL_{s_l}} M_{t}^{Z,u}\le  m_t-z\,\Big|\sigma(\|Z_s\|,s\in[0,s_l])\right)
\end{align*}
Since $Y$ is independent of $\sigma(\|Z_s\|,s\in[0,s_l]$, we continue the above estimation:
\begin{align*}
		&	\bP\left(M_t^Z\le m_t-z\,\big|\sigma(\|Z_s\|,s\in[0,s_l])\right)\\
		\le &\bP(Y\le m_t-m_{t-s_l}-z)+\bP\left(\max_{u\in\cL_{s_l}} M_{t}^{Z,u}\le m_{t-s_l}|\sigma(\|Z_s\|,s\in[0,s_l])\right)\\
		= &\bP\left(Y\le m_t-m_{t-s_l}-z\right)+\left[\bP(M_{t-s_l}^Z\le m_{t-s_l})\right]^{\|Z_{s_l}\|}.
	\end{align*}
	Thus
	\begin{align*}
		&\bP\left(\|Z_{s_{l-1}}\|\le  z^2<\|Z_{s_{l}}\|,M_t^Z\le m_t-z\right)\\
		\le& \bP\left(\|Z_{s_{l-1}}\|\le  z^2<\|Z_{s_{l}}\|)\mP (Y\le m_t-m_{t-s_l}-z \right)+\left[\bP(M^Z_{t-s_l}\le m_{t-s_l})\right]^{z^2}.
	\end{align*}
Since $t-s_l\ge t-s_K\ge t-a^* z\ge t-\frac{(a^*)^2}{\eta}\sqrt{t}\to \infty$
 as $t\to\infty$
	and $\lim_{t\to\infty}\bP(M_t^Z\le m_t)\in(0,1)$, there exist $t(\eta)>1$ and $c_0>0$
	such that for all $t>t(\eta)$,
	\begin{equation}\label{1.2}
		[\bP(M^Z_{t-s_l}\le m_{t-s_l})]^{z^2}\le e^{-c_0z^2}.
	\end{equation}
	As $m_t-m_{t-s_l}-z\le -z(1-\sqrt{2\alpha}\eta l)$ and
$\sqrt{2\alpha}\eta l
 \le\sqrt{2\alpha}\eta K
\le \frac{2\alpha(\rho-1)}{q}
=\frac{2}{\rho+1}<1$, we have by Lemma \ref{Ztlek},
	\begin{align}\label{1.3}
		&\bP\left(\|Z_{s_{l-1}}\|\le  z^2<\|Z_{s_{l}}\|\right)\bP \left(Y\le m_t-m_{t-s_l}-z \right)\nonumber\\
		\le& z^2e^{-q(l-1) \eta z}\mP\left(B_{s_l}\le -z(1-\sqrt{2\alpha}\eta l)\right)\nonumber\\
			=& z^2e^{-q(l-1) \eta z}\mP\left(B_{1}\ge\sqrt{ z}\Big(\frac{1}{\sqrt{\eta l}}-\sqrt{2\alpha}\sqrt{\eta l}\Big)\right)\nonumber\\
		\le &\frac{1}{\sqrt{2\pi}}z^{3/2}e^{-q(l-1)\eta z}
		\Big(\frac{1}{\sqrt{\eta l}}-\sqrt{2\alpha}\sqrt{\eta l}
		\Big)^{-1}
		e^{-\frac{1}{2}(\frac{1}{\sqrt{\eta l}}-\sqrt{2\alpha}\sqrt{\eta l})^2z}\nonumber\\
		\le& \frac{1}{\sqrt{2\pi}}  \left(\frac{1}{\sqrt{a^*}}-\sqrt{2\alpha}\sqrt{a^*}\right)^{-1}z^{3/2}e^{q\eta z}e^{-\sqrt{2\alpha}(\rho-1)z}.
	\end{align}
	Here
	in the second inequality we used
\eqref{est:bm},
and in the final inequality we used the facts that
$\eta l\le \eta K\le a^*$
and $\frac{1}{2}(\frac{1}{\sqrt{\eta l}}-\sqrt{2\alpha}\sqrt{\eta l})^2+q\eta l=(\alpha+q)\eta l+\frac{1}{2l\eta }-\sqrt{2\alpha}\ge \sqrt{2\alpha}(\rho-1)$.

	Combining \eqref{1.4}-\eqref{1.3}, we get that for any $\eta \in(0, a^*/2)$,
	there exist $t_\eta$ and $c_0, C>0$ such that for  $t>t_\eta+t_0$ and  $z\ge 1$,
	$$\bP(M_t^Z\le m_t-z)\le C  (z^2e^{q\eta z}e^{-\sqrt{2\alpha}(\rho-1)z}  + e^{-c_0z^2}).
	$$
	Since $z^2\le 2(q\eta)^{-2}e^{q\eta z}$, and $c_0z^2\ge qa^*z-\frac{(qa^*)^2}{4c_0}$, thus
	$$\bP(M_t^Z\le m_t-z)\le C\left(2(q\eta)^{-2}+e^{\frac{(qa^*)^2}{4c_0}}\right) e^{2q\eta z}e^{-\sqrt{2\alpha}(\rho-1)z} .$$
The proof is now complete.
\hfill$\Box$

\begin{lemma}\label{fact7}
	For any $x\in(0,1)$ and $c\in\R$,
	$$q(1-x)+\frac{\alpha(x-c)^2}{1-x}\ge 2\alpha(\rho-1)(1-c)+\alpha\rho^2\Big(1-\frac{1-c}{\rho}-x\Big)^2.$$
	\end{lemma}
{\bf Proof:}
	Note that the function $(0, \infty)\ni x\to g(x)=a_1^2x+\frac{a_2^2}{x}$ achieves
	its minimum $2a_1a_2$ at the point $x=a_2/a_1$ and for any $x>0$,
	\begin{align}\label{fact6}
		g(x)= 2a_1a_2+\frac{a_1^2}{x}(x-a_2/a_1)^2.
	\end{align}
Then we have that for any $x\in(0,1)$
	\begin{align*}
		q(1-x)+\frac{\alpha(x-c)^2}{1-x}&=(\alpha+q)(1-x)+\frac{\alpha(1-c)^2}{1-x}-2\alpha(1-c)\nonumber\\
		&=\alpha[\rho^2(1-x)+\frac{(1-c)^2}{1-x}-2(1-c)]\nonumber\\
		&=\alpha\Big[2(\rho-1)(1-c)+\frac{\rho^2}{1-x}\left(1-\frac{1-c}{\rho}-x\right)^2\Big]\nonumber\\
		&\ge 2\alpha(\rho-1)(1-c)+\alpha\rho^2\Big(1-\frac{1-c}{\rho}-x\Big)^2,
	\end{align*}
where in the third equality we used \eqref{fact6}.
\hfill$\Box$

	\end{doublespace}

	\vskip 0.2truein
\vskip 0.2truein

\noindent{\bf Yan-Xia Ren:} LMAM School of Mathematical Sciences \& Center for
Statistical Science, Peking
University,  Beijing, 100871, P.R. China. Email: {\texttt
yxren@math.pku.edu.cn}

\smallskip
\noindent {\bf Renming Song:} Department of Mathematics,
University of Illinois,
Urbana, IL 61801, U.S.A.
Email: {\texttt rsong@illinois.edu}

\smallskip

\noindent{\bf Rui Zhang:}  School of Mathematical Sciences \& Academy for Multidisciplinary Studies, Capital Normal
University,  Beijing, 100048, P.R. China. Email: {\texttt
zhangrui27@cnu.edu.cn}
\end{document}